\documentclass[a4paper,12pt]{amsart}

\usepackage{geometry}

\usepackage{amsmath, amsthm, amssymb, mathabx, graphicx}
\usepackage[pdftex]{hyperref}
\usepackage{color}

\newtheorem{theorem}{Theorem}
\newtheorem*{theorem*}{Theorem}

\newtheorem{corollary}{Corollary}
\newtheorem*{corollary*}{Corollary}

\newtheorem*{definition*}{Definition}
\newtheorem{lemma}{Lemma}
\newtheorem*{lemma*}{Lemma}

\newtheorem{proposition}{Proposition}
\newtheorem*{proposition*}{Proposition}

\numberwithin{equation}{section}

\def\supp{{\text{supp }}}
\def\C{{\mathbb C}}
\def\R{{\mathbb R}}

\def\E{{\mathbb E}}

\def\P{{\mathbb P}}

\def\<{{\langle}}
\def\>{{\rangle}}
\def\Var{{\text{Var}}}

\title[Central limit theorems for real zeros of Weyl polynomials]{Central limit theorems for the real zeros of  Weyl polynomials}
\author{Yen Do}
\address{(Yen Do) Department of Mathematics, The University of Virginia, Charlottesville, VA 22904-4137}
\email{yendo@virginia.edu}
\thanks{Y. Do partially supported by NSF grant DMS--1521293.}

\author{Van Vu}
\address{(Van Vu) Department of Mathematics, Yale University, New Haven, CT 06520-8283}
\email{van.vu@yale.edu}
\thanks{V. Vu partially supported by NSF grant  DMS-0901216}
\thanks{This work was initiated  during our visit at the Vietnam Institute for Advanced Study in Mathematics during 2016. We  thank the Institute for the hospitality and support.}

\subjclass[2000]{30B20}

\date{\today}

\begin{document}
\begin{abstract}
We establish the central limit theorem for  the   number of real roots of the  Weyl polynomial $P_n(x)=\xi_0 + \xi_1 x+ \dots + \frac 1{\sqrt {n!}} \xi_n x^n$, where $\xi_i$ are iid Gaussian random variables. The main ingredients in the proof are  new estimates for the correlation functions of the real roots of $P_n$ and a comparison argument exploiting local laws and repulsion properties of these real roots. 
\end{abstract}
\maketitle

\section{Introduction}
 In this paper, we discuss random polynomials with Gaussian coefficients, namely, polynomials of the form 
 
 $$P_n (x)= \sum_{i-0}^n c_i \xi_i x^i $$ where $\xi_i$ are iid standard normal random variables, and $c_i$ are real, deterministic coefficients 
 (which can depend on both $i$ and $n$). 
 
 The central object in the theory of random polynomials, starting with the classical works of Littlewood and Offord \cite{LO1, LO2, LO3}, is the distribution of the real roots. This will  be the focus of our paper.   In what follows, we denote by $N_n$ the number of real roots of $P_n$. 
 
 One important case is when $c_1= \dots = c_n=1$. In this case, the polynomial is often referred to as Kac polynomial. 
  Littlewood and Offord \cite{LO1, LO2, LO3} in the early 1940s, to the surprise of many mathematicians of their time,  showed that $N_{n}$ is typically polylogarithmic in $n$.
 
 \begin{theorem} [Littlewood-Offord]
 	 For Kac polynomials,
 	
 	$$ \frac{\log n} {\log \log n}  \le N_{n} \le \log^2 n$$ with probability $1-o(1)$. 
 	
 \end{theorem}

 Almost simultaneously, Kac \cite{K1} discovered his famous formula for the density function  $\rho(t)$  of $N_{n}$; he show

 \begin{equation} \label{Kac0} \rho(t) =   \int_{- \infty} ^{\infty} |y| p(t,0,y) dy, \end{equation}   where 
 $p(t,x,y)$ is the joint probability density for $P_{n}  (t) =x$ and the derivative $P'_{n} (t) =y$.

 Consequently, 
 \begin{equation} \label{Kacformula}     \E N_{n}  = \int_{-\infty}^{\infty}   dt \int_{- \infty} ^{\infty} |y| p(t,0,y) dy.
 \end{equation}
 
 For Kac polynomials,  he computed $p(t,0,y)$ explicitly and showed \cite{K1} 
 
 \begin{equation}  \label{Kacformula2}   \E N_{n}  = \frac{1}{\pi} \int_{-\infty} ^{\infty} \sqrt { \frac{1}{(t^2-1) ^2} + \frac{(n+1)^2 t^{2n}}{ (t^{2n+2} -1)^2}  }  dt       =   (\frac{2}{\pi} +o(1)) \log n. \end{equation} 
 
More elaborate analysis of  Wilkins \cite{wilkins1988} and also Edelman and Kostlan \cite{ek1995}  
provide a precise estimate of the RHS, showing 

 \begin{equation}  \label{Kacformula3}   \E N_{n}  = \frac{2}{\pi} \log n + C +o(1), \end{equation}  where $C= 0.65...$ is an explicit constant. 
 
  The problem of estimating the variance and establishing the limiting law 
  has  turned out to be significantly harder. Almost thirty years after Kac's work, Maslova solved this problem.
 
 \begin{theorem} \label{theorem:Maslova} \cite{maslova1974, maslova1975}
 	Consider Kac polynomials.  We have, as $n$ tends to infinity
 	
 	\begin{equation*}
 	\frac{N_n-\E N_n}{(Var N_n)^{1/2}} \to N(0,1) .
 	\end{equation*}

 	Furthermore $Var N_n = (K+o(1)) \log n $, where $K= \frac{4}{\pi} (1- \frac{2}{\pi})$. 
 	
 	\end{theorem}

Both Kac's and Maslova's results hold in a more general  setting 
where the gaussian variable is replaced by any random variable with the same mean and variance; see \cite{maslova1974, maslova1975, im1968}.

Beyond the case $c_1 = \dots =c_n=1$, the expectation of $N_n$ is known for many other settings,
see for instance  \cite{bs1986, farahmand1998, maslova1975, ek1995} and the references therein and also the introduction of  \cite{dnv2017} for a recent update.  In many cases, the order of magnitude of  the coefficients $c_i$ (rather than their precise values)
already determines the expectation  $\E N_n$ almost precisely (see the introduction of \cite{dnv2017}).

The limiting law is a more  challenging  problem, and progress has been made only very recently, almost 
40 years after the publication of Maslova's result.  In 2015, Dalmao  \cite{dalmao2015} established the CLT for 
Kostlan-Shub-Smale  polynomials (the case when $c_i = \sqrt { n \choose i} $).
It has been known that in the case the expectation $\E N_n$ is precisely 
$2 \sqrt n$ \cite{ek1995}. 

 \begin{theorem} \label{theorem:Dal} \cite{dalmao2015} 
 	Consider Kostlan-Shub-Smale  polynomials.  We have, as $n$ tends to infinity
 	
 	\begin{equation*}
 	\frac{N_n-\E N_n}{(Var N_n)^{1/2}} \to N(0,1) .
 	\end{equation*}

 	Furthermore $Var N_n = (K+o(1)) \sqrt  n $, where $K= 0.57...$ is an explicit constant. 
 	
 \end{theorem}

There are also many recent  results on random trigonometric polynomial; see  \cite{gw2011, al2013, adl2016, ss2012}; in fact, \cite{dalmao2015} is closely related to \cite{adl2016}, and the proof of Theorem 
\ref{theorem:Dal} used the ideas developed for random trigonometric polynomials from \cite{adl2016}. In particular, 
the papers mentioned above  made essential use of properties of gaussian processes.

In this paper, we  first establish the central limit theorem for $N_n$ for another important  class of random polynomials,  the Weyl polynomials

$$P_n(x) = \sum_{k=0}^n \frac{\xi_k}{\sqrt{k!}} x^k.$$

\begin{theorem}\label{t.CLTall} 
Consider  Weyl polynomials.  We have, as $n \to \infty$,

\begin{equation*}
\frac{N_n-\E N_n}{(Var N_n)^{1/2}} \to N(0,1) .
\end{equation*}

Furthermore $Var N_n = (2K+o(1)) \sqrt n$, where $K= 0.18198..$ is an explicit constant. 
\end{theorem}

It is well known that for Weyl polynomials $\E N_n= (\frac{2}{\pi} +o(1)) \sqrt n $ \cite{ek1995, tv2015}. We give the  exact value of $K$ in the next section. 

Our method for proving the CLT is new, and it  actually yields a stronger result, which establishes 
the following CLT for a very general class of  linear statistics. 

To fix notation, 
let $h:\R\to \R$.  Given $0< \alpha \le 1$, we say that  $h$ is $\alpha$-H\"older continuous on an interval $[a, b]$ if  
$|h(x)-h(y)|\le C|x-y|^\alpha$ for any $a\le x, y \le b$, and the constant $C$ is uniform over $x,y$. 
Below let $Z_n$ denote the (multi)set of the real zeros of $P_n$.

\begin{theorem} \label{t.CLTfinite} There is a finite positive constant $K$ such that the following holds. 
	
	Let $h:\R\to\R$ be bounded, nonzero, and supported on $[-1,1]$ such that
	
	(i) $h$ has finitely many discontinuities and 
	
	(ii) $h$ is  H\"older continuous  when restricted to each interval in the partition of $[-1,1]$ using these discontinuities.
	
	Let  $(R_n) \to \infty$ such that $R_n \le n^{1/2}+o(n^{1/4})$ and let $N_n=\sum_{x\in Z_n} h(x/R_n)$. 
	
	Then 
	$$\lim_{n\to\infty} \frac{\Var [N_n]}{R_n \|h\|_2^2} = K \ \ .$$
	Furthermore, as $n\to\infty$ we have the following convergence in distribution:
	\begin{equation}\label{e.CLTfinite}
	\frac{N_n  - \E N_n}{(\Var \ N_n)^{1/2}} \to N(0,1) \ \ .
	\end{equation}
\end{theorem}

It seems to us that  the approach relying on the properties of  Gaussian processes 
used in  the above mentioned papers  is not applicable in the setting of Theorem  \ref{t.CLTfinite}, with a  general
test function $h$. 

Taking  $h=1_I$ where $I$ is union of finitely many intervals in $[-1,1]$, we obtain the following corollary, which establishes 
the CLT for the number of real roots in  unions of intervals with total length tending to infinity.  

\begin{corollary}\label{c.CLTfinite}
	There is a finite positive constant $K$ such that the following holds. 
	Let $I\subset[-1,1]$ be union of finitely many intervals. Let  $(R_n) \to \infty$ such that $R_n \le n^{1/2}+o(n^{1/4})$ and let $N_n$ be the number of zeros of $P_n$ in $R_n I=\{R_n x, x \in I \}$.
	
	Then 
	$$\lim_{n\to\infty} \frac{\Var [N_n]}{R_n |I|} = K \ \ .$$
	Furthermore, as $n\to\infty$ we have the following convergence in distribution:
	\begin{equation*}
	\frac{N_n  - \E N_n}{(\Var \ N_n)^{1/2}} \to N(0,1) \ \ .
	\end{equation*}
\end{corollary}

For the special case when $I$ is  an interval  of the form $[-a,a]$,  the above asymptotics for the variance of the number of real roots was obtained in \cite{sm2008}.

The  assumption $(R_n)\to \infty$ on the length is optimal, since  asymptotic normality does not hold for intervals of bounded length, due to the repulsion between  nearby real roots.  A similar result of this type was obtained by Granville and Wigman \cite{gw2011} for random trigonometric polynomials, in the special case where the union $I$ consists of one interval.

\section {A sketch of our argument and the outline of the paper }

The heart of the matter is Theorem \ref{t.CLTfinite}.  
It is well-known that most of the real roots of the Weyl polynomial (which we will  denote by $P_n$  in the rest of the proof) 
  are  inside $[-\sqrt n, \sqrt n]$; see for instance \cite{ek1995, tv2015} (see also Lemma~\ref{l.localPn} of the current paper for a local law for the number of real roots of $P_n$). Instead of considering $N_n$, we restrict to the number of real roots  inside $[-\sqrt n, \sqrt n]$. By Theorem \ref{t.CLTfinite}, this variable satisfies CLT. 
   To conclude the proof of Theorem \ref{t.CLTall}, we will  use a tool from 
\cite{tv2015} to bound the number of roots outside this interval, and show that this extra factor is negligible with respect to the 
validity of the CLT.

In order to establish Theorem~\ref{t.CLTfinite}, we first prove a central limit theorem for the random Weyl series
$$P_\infty(x) = \sum_{k=0}^\infty \frac{\xi_k}{\sqrt{k!}} x^k  \ \ .$$
Let $Z$ denote  the (multi)set of the  real zeros of $P_\infty$
where each element in $Z$ is repeated according to its multiplicity. 

For  $h:\R \to \R$ and $R>0$ let $n(R,h) = \sum_{x\in Z} h(x/R)$.
\begin{theorem}\label{t.clt} There is a finite positive constant $K$ such that the following holds.

Let $h:\R\to\R$ be nonzero compactly supported and bounded. Then
$$\lim_{R\to\infty} \frac{\Var[n(R,h)]}{R\|h\|_2^2} = K$$
and as $R\to\infty$ we have the following convergence in distribution:
$$\frac{n(R,h)- \E n(R,h)}{\sqrt{\Var[n(R,h)]}} \to N(0,1)  \ \ .$$
Furthermore, for any $k\ge 1$ it holds that $\E [n(R,h)^k] \le C_{h,k} R^{k}$.
\end{theorem}

The constant $K$ is the same in Theorem~\ref{t.CLTall},   Theorem~\ref{t.CLTfinite}, and Theorem~\ref{t.clt},  and could be computed explicitly: 
\begin{equation}\label{e.Kvalue}
K=\frac 1 \pi + \int (\rho(0,t)-\frac 1 {\pi^2})dt
\end{equation}
where $\rho(s,t)$ is the   two-point correlation function for the real zeros of $P_\infty$. In fact, numerical computation of $K$ was done by Schehr and Majumdar \cite{sm2008} using an explicit evaluation of $\rho(s,t)$ (from the  Kac-Rice formula), giving $K=0.18198...$. For the convenience of the reader and to keep the paper self-contained we sketch some details in Appendix~\ref{s.manju}.

  We deduce Theorem~\ref{t.CLTfinite}   from Theorem~\ref{t.clt} via a comparison argument. Roughly speaking, 
  we try to show that, restricted to certain intervals,  there is a bijection between the real roots of the two functions.
   This argument relies critically on the  repulsion properties of the real roots of $P_n$ and $P_\infty$ (see Section~\ref{s.proveCLTfinite}).    
  
   The rest of the paper is devoted to proving Theorem~\ref{t.clt}. By extending the polynomial to the full series, we 
   can take advantage of the invariance  properties of the root process. The main ingredients of the proof  are  estimates for the correlation functions of the real zeros of $P_\infty$. These correlation function estimates are inspired by related results for the complex zeros of $P_\infty$ by Nazarov and Sodin  \cite{ns2012}, and we adapt  their approach to the real setting. One of the essential steps in \cite{ns2012} is to use a Jacobian formula (which relates the distribution of the coefficients of a polynomial to the distribution of its complex roots)  to estimate the correlation functions of  random polynomials with fixed degrees. Such formula is, however, not available for real roots, and to overcome this difficulty we use  a general expression for correlation functions of real roots of random polynomials due to G\"otze, Kaliada, Zaporozhets \cite{gkz2015}. This expression turns out to be useful to study correlation of small (real) roots, and to remove the smallness assumption we appeal to various invariant properties of the real roots of $P_\infty$.

\vskip2mm 

{\it Acknowledgement.} We would like to thank Manjunath Krishnapur for useful suggestions and Gregory Schehr for a  correction concerning the computation of the explicit constant $K$.


\subsection*{Outline of the paper}
In Section~\ref{s.repulsion} we will prove several  estimates concerning the repulsion properties of the real zeros of $P_n$ and $P_\infty$. In Section~\ref{s.localPn} we will prove some local estimates for the real roots of $P_n$. In Section~\ref{s.proveCLTfinite} we will use these estimates to  prove Theorem~\ref{t.CLTfinite} assuming the validity of Theorem~\ref{t.clt}.

In Section~\ref{s.correlation-functions} we summarize the new estimates for the correlation functions for the real zeros of $P_\infty$, which will be used in Section~\ref{s.proof-clt} to prove Theorem~\ref{t.CLTfinite}.

The proof of the correlation function estimates stated in Section~\ref{s.correlation-functions} will be presented in the remaining sections. 

\subsection*{Notational convention}
By $A\lesssim B$ we mean that there is a finite positive constant $C$ such that $|A|\le CB$. By $A\lesssim_{t_1,t_2\dots,} B$ we mean  that there is a finite positive constant $C$ that may depend on $t_1,t_2,\dots$ such that $|A|\le CB$. Sometimes we also omit the subscripts when the dependency is clear from the context. 

We also say that an event holds with overwhelming probability if it holds with probability at least $1-O_C(n^{-C})$ where $C>0$ is any fixed constant.

For any $I\subset \R$ we'll let $N_n(I)$ be the number of real roots of $P_n$ in $I$.

\section{Real root repulsion} \label{s.repulsion}

In this section we will prove some repulsion estimates for the real roots of $P_n$ (and $P_\infty$). These estimates will be used to  deduce Theorem~\ref{t.CLTfinite} from Theorem~\ref{t.clt}

\subsection{Uniform estimates for $P_n$ and $P_\infty$}
We first establish several basic estimates for the derivatives of $P_n$ and $P_\infty$.  For convenience of notation let $P_{>n} = P_\infty- P_n$.

\begin{lemma} \label{l.derivatives-PnPinfty} Let $I_n = [-n^{1/2}  + n^{1/6}\log n, \ \  n^{1/2} -  n^{1/6}\log n ]$.

For any  $m\ge 0$ integer and $C>0$ there is a constant $c=c(m,C)>0$ such that for any $N>0$ and $n\ge 1$
\begin{eqnarray}
\label{e.Pn-derivative}
\P(\sup_{y\in I_n} |e^{-y^2/2}P^{(m)}_n(y)| > N n^{m/2}) &\lesssim& n e^{-c N^2} \ \ , \\
\label{e.Pinfty-derivative}
\P(\sup_{y\in I_n} |e^{-y^2/2}P^{(m)}_\infty(y)| > N n^{m/2}) &\lesssim& n e^{-cN^2} \ \ , \\
\label{e.P>n-derivative}
\P(\sup_{y\in I_n} |e^{-y^2/2}P^{(m)}_{>n}(y)|  > N n^{-C}) &\lesssim& ne^{-cN^2} \ \ .
\end{eqnarray}
The implicit constants may depend on $C$ and $m$.
\end{lemma}

\proof
Without loss of generality we may assume $N>1$.

We first show  \eqref{e.Pn-derivative}. For any fixed $y$, we have
\begin{eqnarray*}
\Var[e^{-y^2/2}P_n^{(m)}(y)] 
&=& e^{-y^2}\sum_{k=m}^n (k+m)^2\dots(k+1)^2\frac{y^{2k}}{(k+m)!} \\
&=& e^{-y^2} \sum_{k=0}^{n-m} (k+m)\dots(k+1) \frac{y^{2k}}{k!}\\
&<& n^m
\end{eqnarray*}
Since $e^{-y^2/2}P_n^{(m)}(y)$ is centered Gaussian, it follows that for each fixed $y$ we have
$$\P(|e^{-y^2/2}P_n^{(m)}(y)|\ge N n^{m/2})\lesssim e^{-N^2/4}$$

Let $X=(\xi_0,\dots, \xi_n)$ and let $\|.\|$ denote the $\ell_2$ norm on $\R^{n+1}$. By Cauchy-Schwarz, we have the deterministic estimate
$$|e^{-y^2/2} P_n^{(m)}(y)| \le \|X\| \Var[e^{-y^2/2}P_n^{(m)}(y)]  < \|X\| n^{m/2}$$

Let $\delta \in (0,1)$  to be chosen later. Divide the interval $I_n$ into $O(n^{1/2}\delta^{-1})$ intervals of length at most $\delta$. Let $K$ be the collection of the midpoints of these intervals, then by an union bound we have
\begin{eqnarray*}
\P(\exists y\in K: |e^{-y^2/2} P_n^{(m)}(y)| > N n^{m/2})  \le n^{1/2}\delta^{-1} e^{-N^2/4}
\end{eqnarray*}

For any $y'\in I_n$, then there is $y\in K$ such that $|y'-y|\le \delta$. Now, for any $\epsilon \in (0,1)$, using the mean value theorem we have
\begin{eqnarray*}
&&|e^{-(y+\epsilon)^2/2} P_n^{(m)}(y+\epsilon) - e^{-y^2/2} P_n^{(m)}(y)|\\ 
&\le & e^{-(y+\epsilon)^2/2} |P_n^{(m)}(y+\epsilon) -  P_n^{(m)}(y)|  + |[e^{-(y+\epsilon)^2/2}  - e^{-y^2/2}] P_n^{(m)}(y)| \\
&\le&  \epsilon e^{-(y+\epsilon)^2/2}\sup_{\alpha\in (y,y+\epsilon)} |P_n^{(m+1)}(\alpha)|  + \epsilon (y+\epsilon) e^{-y^2/2}  |P_n^{(m)}(y)| \\
&\le& \|X\| \epsilon  [n^{(m+1)/2} + (y+\epsilon)n^{m/2}]\\
&\lesssim&  \epsilon n^{(m+1)/2} \|X\|
\end{eqnarray*}
One could crudely estimate $P(\|X\|> e^{N^2/8}) \le e^{-N^2/4} \E \|X\|^2 = (1+n) e^{-N^2/4}$. (There are sharper estimates for $X$  which follows the chi-squared distribution, but the above estimate is good enough for our purposes.) Therefore by letting $\delta=Nn^{-1/2} e^{-N^2/8}$ and conditioning on the event $\|X\| \le e^{N^2/8}$ we obtain
\begin{eqnarray*}
&&\P(\exists y \in I_n:  |e^{-y^2/2} P_n^{(m)}(y)| > 4N n^{m/2}) \\
&\le&  n N^{-1} e^{-N^2/8} + (n+1)e^{-N^2/4}\\
&\lesssim& n e^{-N^2/8}
\end{eqnarray*}
This completes the proof of \eqref{e.Pn-derivative}.

By the triangle inequality it remains to show \eqref{e.P>n-derivative}. 

We proceed as before. Given any fixed $y$   we have
\begin{eqnarray*}
\Var[e^{-y^2/2}P_{>n}^{(m)}(y)] 
&=&   e^{-y^2}\sum_{k=n-m+1}^\infty (k+1)^2\dots (k+m)^2\frac{y^{2k}}{(k+m)!}\\
&\lesssim&  y^{2m} e^{-y^2} \sum_{k=n-m+1}^\infty \frac{y^{2(k-m)}}{(k-m)!}\\
\end{eqnarray*}

Let $y_0=n^{1/2}-n^{1/6} \log n$. Then for $n$ large enough (relative to $m$) we have $n-m+1\ge y_0^2> y^2$, consequently  for each $k\ge n-m+1$ the function $h(y)=2k\log |y| - y^2$ is increasing over $y\in (0,y_0$.  It follows that
\begin{eqnarray*}
\Var[e^{-y^2/2}P_{>n}^{(m)}(y)]  &\lesssim& n^{m} e^{-y_0^2}\sum_{k=n-2m+1}^\infty \frac{y_0^{2k}}{k!}
\end{eqnarray*}
Since $\sqrt n\lesssim y_0\le \sqrt n$ and $m=(1)$, it follows that
\begin{eqnarray*}
\Var[e^{-y^2/2}P_{>n}^{(m)}(y)]  &\lesssim& n^{m} e^{-y_0^2}\sum_{k=n}^\infty \frac{y_0^{2k}}{k!}\\
&\lesssim& n^{m} e^{-y_0^2}\sum_{k=n}^{[(1.01)n]} \frac{y_0^{2k}}{k!} \\
&\lesssim& n^{m+1} e^{-y_0^2} \frac{y_0^{2n}}{n!} 
\end{eqnarray*}
(here we used the fact that $y_0^2/k < 1/1.01<1$ if $k\ge (1.01)n$, and $y^2/k\le 1$ for $k\ge n$). Consequently,
\begin{eqnarray*}
\Var[e^{-y^2/2}P_{>n}^{(m)}(y)]  &\lesssim& n^{m+1} e^{-y_0^2} \frac{y_0^{2n}}{n!} \\
&\lesssim& n^{m+1} e^{-y_0^2} \frac{y_0^{2n}}{(n/e)^{n+1/2}} \\
&=& n^{m+1/2}  e^{-y_0^2+ 2n \log y_0 - n\log n + n}
\end{eqnarray*}
Now, for brevity write $y_0=\sqrt n (1-\beta)$ where $\beta=n^{-1/3}\log n = o(1)$, then
\begin{eqnarray*}
&& -y_0^2+ 2n \log y_0 - n\log n + n \\
&=& -n(1-2\beta+\beta^2)  + n \log n +2n\log(1-\beta) - n\log n +n \\
&=& 2n[\beta - \frac{\beta^2}2 + \log (1-\beta)]\\
&=& -2n \frac{\beta^3}3 (1+O(\beta)) \\
&\le& -n \beta^3/3 = -\log n^3/3
\end{eqnarray*}
when $n$ is large. Therefore
\begin{eqnarray}
\label{e.variance-P>n}
\Var[e^{-y^2/2}P_{>n}^{(m)}(y)]   &\lesssim&  n^{m+1/2}  e^{-\log n^3/3} \lesssim_{C,m} n^{-C}
\end{eqnarray}
for any $C>0$. Therefore for any fixed $y$ such that $|y|\le \sqrt n(1-\log n/n^{1/3})$ we have
\begin{eqnarray*}
\P(|e^{-y^2/2}P_{>n}^{(m)}(y)|> N n^{-C}) \lesssim e^{-N^2/4}
\end{eqnarray*}

Let $X=(\xi_{n+1},\xi_{n+2},\dots, \xi_{3n}, \xi_{3n+1}/2, \dots, \xi_m/2^{m-3n},\dots)$ and let $\|.\|$ denote the $\ell_2$ norm on $\R^{n+1}$. By Cauchy-Schwarz and using $y^2\le n$, we have the deterministic estimate
\begin{eqnarray*}
|e^{-y^2/2} P_{>n}^{(m)}(y)|
&\le& \|X\| e^{-y^2/2} (\sum_{k=n+1-m}^{3n-m} \frac{(k+1)\dots (k+m)y^{2k}}{k!} + \\
&& \qquad \qquad + \sum_{k>3n-m} \frac{4^{k+m-3n}(k+1)\dots (k+m)y^{2k}}{k!})^{1/2} \\
&\lesssim& \|X\| e^{-y^2/2}(\sum_{k=n+1}^{3n} \frac{k^m y^{2k}}{k!} + O(\frac{n^m y^{6n}}{(3n)!}))^{1/2}\\
&\lesssim& \|X\| n^{m/2} 
\end{eqnarray*}

Let $\delta \in (0,1)$  to be chosen later. Divide the interval $I_n$ into $O(n^{1/2}\delta^{-1})$ intervals of length at most $\delta$. Let $K$ be the collection of the midpoints of these intervals, then by an union bound we have
\begin{eqnarray*}
\P(\exists y\in K: |e^{-y^2/2} P_{>n}^{(m)}(y)| > N n^{-C})  \le n^{1/2}\delta^{-1} e^{-N^2/4}
\end{eqnarray*}

For any $y'\in I_n$, then there is $y\in K$ such that $|y'-y|\le \delta$. Now, for any $\epsilon \in (0,1)$, using the mean value theorem we have

\begin{eqnarray*}
&&|e^{-(y+\epsilon)^2/2} P_{>n}^{(m)}(y+\epsilon) - e^{-y^2/2} P_{>n}^{(m)}(y)|\\ 
&\le & e^{-(y+\epsilon)^2/2} |P_{>n}^{(m)}(y+\epsilon) -  P_{>n}^{(m)}(y)|  + |[e^{-(y+\epsilon)^2/2}  - e^{-y^2/2}] P_{>n}^{(m)}(y)| \\
&\le&  \epsilon e^{-(y+\epsilon)^2/2}\sup_{\alpha\in (y,y+\epsilon)} |P_{>n}^{(m+1)}(\alpha)|  + \epsilon (y+\epsilon) e^{-y^2/2}  |P_{>n}^{(m)}(y)| \\
&\lesssim& \|X\| \epsilon  [n^{(m+1)/2} + (y+\epsilon)n^{m/2}]\\
&\lesssim&  \epsilon n^{(m+1)/2} \|X\|
\end{eqnarray*}
One could crudely estimate $P(\|X\|> e^{N^2/8}) \le e^{-N^2/4} \E \|X\|^2  \lesssim n e^{-N^2/4}$. Therefore by letting $\delta=Nn^{-1/2} e^{-N^2/8}$ and conditioning on the event $\|X\| \le e^{N^2/8}$ we obtain
\begin{eqnarray*}
&&\P(\exists y \in I_n:  |e^{-y^2/2} P_{>n}^{(m)}(y)| \gtrsim N n^{m/2}) \\
&\lesssim&  n N^{-1} e^{-N^2/8} + n e^{-N^2/4}\\
&\lesssim& n e^{-N^2/8}
\end{eqnarray*}
This completes the proof of \eqref{e.P>n-derivative}.

\endproof

\subsection{Replusion of the real roots}

In this section we prove estimates concerning the separation  of real roots of $P_n$ and $P_\infty$ in $I_n =[-n^{1/2}+n^{1/6}\log n, n^{1/2}-n^{1/6}\log n]$.

\begin{lemma}\label{l.drPn} For any $c_2>0$ the following estimates hold for $c_1>c_2+2$:

(i) $ \P\Big(\exists \ x\in I_n: \ \ P_n(x)=0, \ \ |\frac{d}{dx}(e^{-x^2/2}P_n(x))|< n^{-c_1} \Big) \lesssim n^{-c_2}$

(ii) $ \P\Big(\exists \ x, x' \in I_n: \ \   P_n(x)= P_n(x')=0, \ \  0< |x-x'|< n^{-c_1} \Big) \lesssim n^{-c_2}$
\end{lemma}


\proof 

For convenience of notation let $q_n(x) =e^{-x^2/2} P_n(x)$. Clearly $q_n$ and $P_n$ have the same real roots. Furthermore, for $x\in I_n$ it holds that
\begin{eqnarray*}
q_n'(x) &=& e^{-x^2/2}P'_n(x) + (-x)e^{-x^2/2}P_n(x) \\
&\le& e^{-x^2/2}|P_n'(x)| + \sqrt n e^{-x^2/2} |P_n(x)|
\end{eqnarray*}
and similarly
\begin{eqnarray*}
q_n''(x) &=& e^{-x^2/2}P''_n(x) + 2(-x)e^{-x^2/2}P'_n(x)  + (x^2-1)e^{-x^2/2}P_n(x)\\
&\lesssim& e^{-x^2/2}|P''_n(x)| + \sqrt n e^{-x^2/2}|P'_n(x)|  +  n e^{-x^2/2}|P_n(x)|
\end{eqnarray*}
Thus, using Lemma~\ref{l.derivatives-PnPinfty} with $N=  C\log^{1/2} n$ with $C>0$ large, we obtain
\begin{equation}\label{e.P_n''}
\P(\sup_{y\in I_n} |q''_n(y)| \gtrsim n \log^{1/2} n) \lesssim n e^{-cC\log n} < n^{-c_2}  \ \ .
\end{equation}

(i) Let $\delta = n^{-c_1}$. 

Suppose that $q_n(x)=0$ and $|q'_n(x)|<\delta$ for some fixed $x\in I_n$.  Then for every $x'\in I_n$ with $|x'-x|\le \delta$,   conditioning on the event $\sup_{|y|\in I_n} |q_n''(y)| \lesssim n\log^{1/2} n$ and using the mean value theorem, we have 
\begin{eqnarray*}
q_n(x') &=& q_n(x) + (x'-x) q'_n(x)+ O(\frac{(x-x')^2 n\log^{1/2}n}2) \\
&\lesssim& \delta^2 + \delta^2 n \log^{1/2}n \\
&\lesssim& \delta^2 n \log^{1/2}n =:\beta
\end{eqnarray*}
Now, divide the interval $I_n$ into $O(n^{1/2}\delta^{-1})$ intervals of length at most $\delta/2$. Using the above estimates and using an union bound, it follows that
\begin{eqnarray*}
&&\P\Big(\exists \ x\in I_n: \ \ q_n(x)=0, \ \ |q_n'(x)|< \delta \Big) \\
&\lesssim& \sqrt n \delta^{-1} \sup_{x\in I_n} \P(|q_n(x)| \lesssim \beta) + n^{-c_2}
\end{eqnarray*}

Now, for each $x\in I_n$  there is $0\le j \le n$ depending on $x$ such that $|e^{-x^2/2}\frac{x^j}{\sqrt {j!}}| \gtrsim n^{-1/2}$.
To see this, we invoke \eqref{e.variance-P>n} for $m=0$ and obtain
\begin{eqnarray*}
e^{-x^2}\sum_{j>n} \frac{x^{2j}}{j!}  = \Var[e^{-x^2/2}P_{>n}(x)]  \lesssim n^{1/2} e^{-\log^3 n/3}
\end{eqnarray*}
Consequently for $x\in I_n$ we have $\sum_{j=0}^n \frac{x^{2j}}{j!} \gtrsim e^{x^2}$
and therefore one could select a $j\in \overline{0,n}$ with the stated properties.

Given such a $j$, we condition on $e^{-x^2/2} \sum_{i\ne j} \xi_i \frac{x^i}{\sqrt{i!}}$, which is independent from $\xi_j$, obtaining
$$P(|q_n(x)| < \delta^{1+\epsilon}) \le \sup_{z} P(\frac{x^{j}e^{-x^2/2}}{\sqrt{j!}} \xi_j \in (z-\beta, z+\beta))$$ 
$$\lesssim  \frac{e^{x^2/2}\sqrt{j!}}{|x^j|} \beta \lesssim n^{1/2} \beta$$
since the density of the Gaussian distribution (of $\xi_j$) is bounded. Note that the implicit constants are independent of $x\in  I_n$. Consequently,
\begin{eqnarray*}
\P\Big(\exists \ x\in I_n: \ \ P_n(x)=0, \ \ |P_n'(x)|< \delta \Big) 
&\lesssim& n^{-c_2} +  n \delta^{-1}  \beta\\
&=& n^{-c_2} + \delta n^2 \log^{1/2}n\\
&=& n^{-c_2} + n^{2-c_1}\log^{1/2} n \lesssim  n^{-c_2}
\end{eqnarray*}
provided that $c_1>c_2+2$.

(ii)   Assume that for some $x\ne x'$ in $I_n$  we have $q_n(x)=q_n(x')=0$. By the mean value theorem there is some $x''$ between $x, x'$ such that $q'_n(x'')=0$. Let $\delta=n^{-c_1}$ as before. Conditioning on the event $\sup_{|y|\in I_n} |q_n''(y)| \lesssim n\log^{1/2} n$ (which holds with probability $1-O(n^{-c_2})$) and  using the mean value theorem we have
$$q'_n(x) = q'_n(x'') + |x-x''| O(n\log^{1/2}n) = O(\delta n\log^{1/2}n)$$
therefore for any $y\in [x-\delta, x+\delta]$ it holds that
\begin{eqnarray*}
q_n(y) &=& q_n(x)+(y-x)q'_n(x) + (y-x)^2  O(n\log^{1/2}n) \\
&=& O(\delta^2 n\log ^{1/2}n )
\end{eqnarray*}

The rest of the proof similar to (i).  
\endproof

Using an entirely  similar argument, we also have the following series analogue of Lemma~\ref{l.drPn}.

\begin{lemma}\label{l.drPinfty} For any $c_2>0$ the following estimates hold for $c_1>c_2+2$:

(i) $ \P\Big(\exists \ x\in I_n: \ \ P_\infty(x)=0, \ \ |\frac{d}{dx}(e^{-x^2/2}P_\infty(x))|< n^{-c_1} \Big) \lesssim n^{-c_2}$

(ii) $ \P\Big(\exists \ x, x' \in I_n: \ \   P_\infty(x)= P_\infty(x')=0, \ \  0< |x-x'|< n^{-c_1} \Big) \lesssim n^{-c_2}$
\end{lemma}

\section{Local law for $P_n$}\label{s.localPn}
In this section we prove a local law for $P_n$, which will be used in the proof of Theorem~\ref{t.CLTfinite} and Theorem~\ref{t.CLTall}.

\begin{lemma}\label{l.localPn}  The following holds with overwhelming probability: for any interval $I\subset \R$   it holds that
$$N_n(I) \lesssim (1+ |I \cap [-\sqrt n,\sqrt n]|) n^{o(1)}  \ \ .$$ 
\end{lemma}

A variant of  Lemma~\ref{l.localPn} for complex zeros of (non-Gaussian) Weyl polynomials was considered in \cite{tv2015} (see estimates (87,88) of \cite{tv2015}). The proof given below for Lemma~\ref{l.localPn} is inspired by the (complex) argument in \cite{tv2015}. Our setting is simpler because $P_n$ is Gaussian thus our condition on $I$ is weaker (in comparison to the requirement that $I\subset \{n^{-C} \le |x| \le C\sqrt n\}$ in \cite{tv2015}).

We will need the following estimate \cite[Proposition 4.1, arXiv version]{tv2015}; to keep the current paper self-contained we will include a proof of this estimate shortly.

\begin{proposition} \label{p.prop41tv2014} Let $n\ge 1$ be integer and $f$ be a random polynomial of degree at most $n$. Let $z_0\in \C$ be depending on $n$, and let $n^{-O(1)} \lesssim c\le r \lesssim n^{O(1)}$ be quantities that may depend on $n$. 

Let $G:\C \to \C$ be a deterministic smooth function that may depend on $n$ such that
$$\sup_{z\in  B(z_0, r+c) \setminus B(z_0, r-c)} |G(z)| \lesssim n^{O(1)}$$

Assume that for any $z\in B(z_0, r+c) \setminus B(z_0, r-c)$ one has
$$\log |f(z)| = G(z)+O(n^{o(1)})$$ with overwhelming  probability.

Then with overwhelming probability the following holds: $f \not\equiv 0$ and the number $N$ of roots of $f$ in $B(z_0,r)$ satisfies
$$N =\frac 1 {2\pi} \int_{B(z_0,r)} |\Delta G(z)|dz + O(n^{o(1)} c^{-1}r)  + O(\int_{ B(z_0, r+c) \setminus B(z_0, r-c)} |\Delta G(z)| dz)$$
\end{proposition}

\proof We largely follow \cite{tv2015}. We first prove the upper bound
$$N \le \frac 1 {2\pi} \int_{B(z_0,r)} |\Delta G(z)|dz + O(n^{o(1)} c^{-1}r)  + O(\int_{ B(z_0, r+c) \setminus B(z_0, r-c)} |\Delta G(z)| dz)$$

Let $Z$ denote the (multi)set of zeros of $f$. We then estimate
$$N\le \sum_{z\in Z} \phi(z) =\frac 1 {2\pi} \int_{\C} \Delta \phi(z) \log |f(z)|dz$$
$$=\frac 1 {2\pi}\int_{\C} \Delta \phi(z) G(z)dz + O(\int_{\C} |\Delta \phi(z)|  |h(z)|dz)$$
$$=\frac 1 {2\pi}\int_{\C}  \phi(z) \Delta G(z)dz + O(\int_{\C} |\Delta \phi(z)|  |h(z)|dz)$$
Now it is clear that
$$\int_{\C}   \phi(z)  \Delta G(z)dz = \int_{B(z_0,r)} \Delta G(z) dz + O(\int_{ B(z_0, r+c) \setminus B(z_0, r-c)} |\Delta G(z)| dz)$$
therefore to estabish the upper bound it remains to show that
$$\int_{B(z_0,r+c)\setminus B(z_0,r-c)}   |h(z)|dz \lesssim n^{o(1)} cr$$
with overwhelming probability.

We first show that 
\begin{eqnarray}\label{e.2mMC}
\int_{B(z_0,r+c)\setminus B(z_0,r-c)} |h(z)|^2 dz=O(n^{O(1)})
\end{eqnarray}
with overwhelming probability. Since $h=\log |f|-G$ and $|G(z)|=O(n^{O(1)})$ uniformly on $B(z_0,r+c)\setminus B(z_0,r-c)$ and $r,c$ are polynomial in $n$,  it remains to show the same estimate for $\log|f|$.

Now, let $z_1=z_0+r$. By the given hypothesis, with overwhelming probability $\log |f(z_1)| =G(z_1)+O(n^{o(1)}) = O(n^{O(1)})$ (which in particular means that $z_1$ is not one of the zeros of $f$ on this event). Conditioning on this event, it remains to show that
$$\int_{B(z_0,r+c)\setminus B(z_0,r-c)} |\log |f(z)| - \log |f(z_1)||^2 dz =O(n^{O(1)})$$
Factorizing $P_n(z)=C\prod_{k=1}^n(z-\xi_k)$ where $(\xi_k)$ are the complex roots of $f$, it suffices to show that for any $\alpha\in \C\setminus z_1$ we have
$$\int_{B(z_0,r+c)\setminus B(z_0,r-c)}  \log^2 |\frac{z-\alpha}{z_1-\alpha}| dz =O(n^{O(1)})$$
(uniformly over $\alpha$). There are two cases: First, if $|\alpha -z_0|>r+2c$ then clearly 
$$|\frac{|z-\alpha|}{|z_1-\alpha|} -1| \le \frac{|z_1-z|}{c} =O(n^{O(1)})$$
uniformly over $z\in B(z_0,r+c)\setminus B(z_0,r-c)$, which implies the desired estimate. Secondly, if $|\alpha-z_0|\le r+2c$ then we estimate $\log^2 |\frac{z-\alpha}{z_1-\alpha}|  \lesssim \log^2|z-\alpha| + \log^2 |z_1-\alpha|$ and the desired estimate follows from 
$$\int_{B(0,2r+3c)} \log^2 |z| dz = O((1+2r+3c)^{O(1)}) = O(n^{O(1)})$$

Now we condition on the event that \eqref{e.2mMC} holds. For $m\ge 1$ let $x_1,\dots, x_m$ be independently randomly selected from $B(z_0,r+c)\setminus B(z_0,r-c)$ (which has measure $4rc$). (These points are also chosen independent of $f$.) Using the Monte Carlo sampling lemma and \eqref{e.2mMC}, we have
$$ \int_{B(z_0,r+c)\setminus B(z_0,r-c)} |h(z)| \le 4rc[\frac 1 m \sum_{k=1}^m |h(x_k)|]+ (m\delta)^{-1/2} n^{O(1)}$$
with probability at least $1-\delta$.  Given any $C>0$ let $\delta=n^{-C}$ and $m=n^{A+C}$ where $A$ is very large. By the given hypothesis, for each $k=1,\dots, m$ we have $h(x_k)| = O(n^{o(1)})$ with probability at least $1-O_{A,C}(n^{-(2A+C)})$. Therefore using union bound with probability at least $1-O(n^{-C})$ we have (still conditioning on \eqref{e.2mMC} holding):
$$ \int_{B(z_0,r+c)\setminus B(z_0,r-c)} |h(z)|  \lesssim rc n^{o(1)} + n^{-A/2}n^{O(1)}$$
Since $A$ could be chosen arbitrarily large and $r,c$ are at least some negative power of $n$, it follows that  with probability at least $1-O(n^{-C})$ we have (still conditioning on \eqref{e.2mMC} holding):
$$ \int_{B(z_0,r+c)\setminus B(z_0,r-c)} |h(z)|  \lesssim rc n^{o(1)}$$
Since \eqref{e.2mMC} holds overwhelmingly, by removing the conditioning it follows that the following estimate hold with probability $1-O(n^{-C})$ (any $C>0$):
$$ \int_{B(z_0,r+c)\setminus B(z_0,r-c)} |h(z)|  \lesssim rc n^{o(1)}$$
in other words it holds overwhelmingly, as desired.

For the lower bound for $N$ we'll choose  $\phi$ to be supported on $B(z_0,r)$ and equal 1  on $B(z_0,r-c)$, the rest of the argument is entirely similar.
\endproof

We'll use a crude estimate for the roots of $P_n$:
\begin{lemma}\label{l.rootboundPn} Given any $C>0$, with probability at least $1-O(n^{-C})$ the   roots of $P_n$ satisfy $|z|\le n^{(3C+2)/ 2} $.
\end{lemma}
\proof  Without loss of generality assume $n\ge 2$. Let $X=(\xi_0,\dots, \xi_{n-1})$ and let $\|.\|$ denote the $\ell_2$ norm on $\R^n$. By Cauchy-Schwartz, we have the deterministic estimate
$$|\sum_{j=0}^{n-1}\xi_j \frac{z^j}{\sqrt {j!}}| \le \|X\| (\sum_{j=0}^{n-1} \frac{|z|^{2j}}{j!})^{1/2}$$
For any $|z|>n^{(3C+ 2)/ 2}$ it is clear that the sequence $(|z|^{2j}/j!)_{j=0}^n$ is lacunary 
$$\frac{|z|^{2j}/j!}{|z|^{2j-2}/(j-1)!} = \frac{|z|^2}{j} \ge n^{3C+2}/j \ge n^{3C+1} >1$$
therefore we have the deterministic bound 
$$(\sum_{j=0}^{n-1}\frac{|z|^{2j}}{j!})^{1/2} \lesssim_C n^{-(3C+1)/2} \frac{|z|^{n}}{\sqrt{n!}} $$
Consequently it suffices to show that  the event $\{\|X\| \le \frac 1 {M} n^{(3C+1)/2}|\xi_n|\}$ has probability at least $1-O_{M,C}(n^{-C})$, any $M>0$. 
Since $\E \|X\|^2 = n$, it follows that 
$$P(\|X\|< n^{(C+1)/2}) = 1-O(n^{-C})$$
and using boundedness of the density of Gaussian we have
$$P(|\xi_n| \ge M n^{-C}) = 1- O_M(n^{-C})$$
thus taking the intersection of these two events we obtain the desired claim.
\endproof

\subsection{Proof of Lemma~\ref{l.localPn}}
We now begin the proof of Lemma~\ref{l.localPn}.  Note that $P_n/|VarP_n|^{1/2}$ is normalized Gaussian. It follows that for any $z$
$$\log |P_n(z)| = \frac 1 2\log |Var P_n(z)| + O(n^{o(1)})$$
with overwhelming probability (the implicit constant is independent of $z$ but the bad event may depend on $z$). And 
$$Var P_n(z) = \sum_{j=0}^n \frac{|z|^{2j}}{j!}$$

Let $z$ be such that $|z|\ge \sqrt n$.  Then the sequence $1 \le |z|^2/1! \le \dots \le |z|^{2n}/n!$ is increasing. It follows that $|z|^{2n}/n! \le Var[P_n] \le (n+1)|z|^{2n}/n!$, and consequently using Stirling's formula we have the uniform bound
$$\log |Var P_n(z)|  = 2n \log |z| - \log (n!) + O(\log n)$$ 
$$ = 2n\log |z| - (n\log n -n) + O(n^{(o(1)})$$

If $|z|\le n^{1/2}$ then $|z|^{2k}/(k)! \ge |z|^{2k+2}/(k+1)!$ for any $k\ge n$ and when $k>2n$ we have $|z|^{2k}/k! \ge 2|z|^{2k+2}/(k+1)!$. Thus $n^{-1} e^{|z|^2} \lesssim Var [P_n] \le e^{|z|^2}$ therefore
$$\log |Var P_n(z)|   =|z|^2 +  O(n^{o(1)})$$

We now take $G(z)=\frac 1 2g(|z|)|z|^2+(1-g(|z|))[n\log |z|  - \frac 1 2n(\log n-1)]$ which is smooth where $g:\R \to [0,1]$ a bump function such that  $g(x)=1$ for $|x|\le \sqrt n$ and $g(x)=0$ for $|x|\ge \sqrt n + 1$.  In the transitional region $\sqrt n \le |z| \le \sqrt n +1$, by examination we have
$$2n\log |z|  - n\log n + n = |z|^2 + O(1)$$
Therefore for each $z$ with overwhelming probability it holds that
$$\log |P_n(z)|  = \frac 1 2 \log |Var P_n(z)| + O(n^{o(1)}) =  G(z)+O(n^{o(1)})$$

Note that $G$ is depending only on $|z|$ and satisfies polynomial bound $G(z)=O(n^{O(1)})$ if $|z|$ is also at most polynomial in $n$. Furthermore, 
$$\Delta G(z) = \begin{cases} 2, & |z|\le \sqrt n\\
0, & |z|\ge \sqrt n +1
\end{cases}$$
and for $\sqrt n < |z| < \sqrt n+1$ using the polar coordinate form  of $\Delta$ it holds that
$$\Delta G(z) = [\frac 1{r} \partial_r +\partial_{rr}] \Big(g(r)\frac {r^2}2+ (1-g(r))(n\log r - \frac 1 2 n \log n + \frac n2)\Big) $$
$$=O(\frac 1 r |g'(r)|) + O(|g''(r)|) + O(|g'(r)| |r-\frac n r|) = O(1)$$

Now, let $C>0$, then by Lemma~\ref{l.rootboundPn} with probability $1-O(n^{-C})$ the roots of $P_n$ satisfy $|z|\le N:=n^{(3C+2)/2}$.  

We now apply Proposition~\ref{p.prop41tv2014} with $z_0=N$, $r=N/2$, and $c=N/4$. Then with overwhelming probability 
$$N_n[N/2, 3N/2] \lesssim  \int_{B(N,N/2)} 1_{B(0,\sqrt n+1)} + O(n^{o(1)}) + \int_{B(N,3N/2)\setminus B(N,N/2)} 1_{B(0,\sqrt n+1)}$$
$$ = O(n^{o(1)})$$
We then repeat (variance of) this argument $O(\log N)$ times with a decreasing lacunary sequence of $z_0$ (starting from $N$). Then  with overwhelming probability, in $[\sqrt n+2, N]$ there are $O(n^{o(1)}\log N) =O(n^{o(1)})$ real roots. By a similar argument, we have the same bound in $[-N, -\sqrt n - 2]$ with overwhelming probability.

We now consider the real roots in $[-\sqrt n-2, \sqrt n+2]$.

Let $z_0=\sqrt n$ and $r=c=2$, it follows that with overwhelming probability the number of real roots in $[\sqrt n-2,\sqrt n+2]$ is $O(n^{o(1)})$. By repeating this argument it follows that for any interval $I_0 \subset [-\sqrt n- 2, \sqrt n+2]$ of length  $1$  with overwhelming probability the number of real roots in $I$ is $O(n^{o(1)})$. Of course if $I_0$ has length less than $1$ then using monotonicity of $N_n(I)$ we also have the same upper bound. Dividing $ [-\sqrt n- 2, \sqrt n+2]$ into intervals of length $1$ and taking the union bound, it follows that one could could ensure that for all subintervals of length $1$ with overwhelming probability.

Consequently, given any $C>0$, with probability $1-O(n^{-C)}$, for any interval $I\subset\R$   we have $$N_n(I) \lesssim (1+ |I \cap [-\sqrt n, \sqrt n]|) n^{o(1)} \ \ .$$
This completes the proof of Lemma~\ref{l.localPn}.

\section{Proof of Theorem~\ref{t.CLTall} assuming Theorem~\ref{t.CLTfinite}}\label{s.proveCLTall}
Recall the notation that   $N_n(I)$ denotes the number of real roots of $P_n$ in $I\subset \R$. Let $h=1_{[-1,1]}$  and $R_n=\sqrt n$. Let $N_{n,in}:=N_n([-\sqrt n,\sqrt n])$ and $N_{n,out}=N_n-N_{n,in}$. Then by Theorem~\ref{t.CLTfinite}, we have
$$Var[N_{n,in}]/2\sqrt n \to K \in (0,\infty)$$
$$\frac{N_{n,in} - \E N_{n,in}}{\sqrt{Var[N_{n,in}]}} \to N(0,1)$$
as $n\to \infty$, and the second convergence is in distribution. By Lemma~\ref{l.localPn}, with overwhelming probability we have $N_{n,out}= O(n^{o(1)})$, and we always have $N_{n,out}\le n$ deterministically. Consequently 
$$ \E N_{n,out}^2 = O(n^{o(1)}) = o(\E N_{n,in}^2)$$
and therefore $Var[N_{n,out}]=O(n^{o(1)})$ and so 
$$Var[N_n] = Var[N_{n,in}] (1+o(1)) = 2\sqrt n K (1+o(1))$$
Furthermore,   with overwhelming probability we have
$$\frac{N_n - \E N_n}{\sqrt{Var[N_n]}} =  o(1)+ \frac{N_{n,in} - \E N_{n,in}}{\sqrt {Var[N_n]}}$$
$$= o(1)  + (N_{n,in} - \E N_{n,in}) \Big[\frac{1}{\sqrt {Var[N_{n,in}]}}+ \frac{O(n^{o(1)})}{\sqrt {Var[N_{n,in}]}   \sqrt {Var[N_{n}]}}\Big]$$
$$=o(1)+  \frac{N_{n,in} - \E N_{n,in}}{\sqrt {Var[N_{n,in}]}}(1+o(1))$$
Thus by Slutsky's theorem (see e.g. \cite[Chapter 7]{ash2000}) it follows that $\frac{N_n - \E N_n}{\sqrt{Var[N_n]}}\to N(0,1)$ in distribution. 

\section{Proof of Theorem~\ref{t.CLTfinite} assuming Theorem~\ref{t.clt}}\label{s.proveCLTfinite}

The comparison argument  in this section is inspired by similar arguments in \cite{dnv2015,nnv2015}.

Recall that $N_n=\sum_{x\in Z_n} h(x/R_n)$ and $Z_n$ is the multiset of the real zeros of $P_n$.

Denote $N_\infty:=n(R,h)=\sum_{x\in Z} h(x/R_n)$ where $Z$ is the multiset of  the real zeros of $P_\infty$.  
Let $N_G=N(0,1)$ be  the standard Gaussian random variable, and
$$N^*_n:=\frac{N_n - \E N_n}{(\Var \ N_n )^{1/2}} \ \ , \ \ N_\infty^* =\frac{N_\infty  - \E N_\infty}{(\Var N_\infty)^{1/2}} \ \ .$$
Applying Theorem~\ref{t.clt}, we obtain
$$\lim_{n\to\infty} \frac{\Var[N_\infty]}{\|h\|_2^2 R_n} = K$$
and $N^*_\infty \to N_G$ in distribution.

To deduce Theorem~\ref{t.CLTfinite}, we will compare $N_n$ with $N_\infty$.
\begin{lemma}\label{l.compareN} As $n\to\infty$, it holds that
$$\E|N_n-N_\infty|^2  = o(R_n)$$
\end{lemma}
Below we prove Theorem~\ref{t.CLTfinite} assuming Lemma~\ref{l.compareN}.

\proof[Proof of Theorem~\ref{t.CLTfinite}]
For convenience let $\Delta N_n = N_n - N_\infty$. 
It follows from Lemma~\ref{l.compareN} that $|\E \Delta_n| = |\E N_n  - \E N_\infty|=o(R_n^{1/2})$. Using the $L^2$ triangle inequality 
we also obtain $|\sqrt{\Var [N_n]}-\sqrt{\Var [N_\infty]}|  = o(R_n^{1/2})$.
Since $\Var[N_\infty]=2R_nK(1+o(1)) = \Theta(\sqrt n)$ as $n\to\infty$, we obtain $\Var[N_n]=2R_nK(1+o(1))$, in particular  $\Var[N_n]  = \Theta(R_n)$.

Now,
\begin{eqnarray*}
N^*_n &=&   \frac{\Delta N_n - \E \Delta N_n}{[Var N_n]^{1/2}} +  \frac {N_\infty- \E N_\infty}{[Var N_n]^{1/2}} \\
&=& O(\frac{|\Delta N_n - \E \Delta N_n|}{R_n^{1/2}}) + N_\infty^*\Big(1+  \frac{\sqrt{Var N_\infty}  - \sqrt{Var N_n}}{ \sqrt{Var N_n]}} \Big)\\
&=& O(\frac{|\Delta N_n|}{R_n^{1/2}} +o(1)) + N_\infty^*\Big(1+o(1)\Big)
\end{eqnarray*}
Since $\E|\Delta_n|^2 = o(R_n)$, it follows that $\frac{|\Delta N_n|}{R_n^{1/2}}\to 0$ in probability. Therefore by Slutsky's theorem (see e.g. \cite[Chapter 7]{ash2000}) it follows that $N^*_n$ conveges to $N(0,1)$ in distribution. 

\endproof

Our proof of Lemma~\ref{l.compareN} will use a comparison argument. 
More specifically, we'll show that with high probability $\sup_{x\in I_n} |P_n(x)-P_\infty(x)|$  is very small in comparison to the typical distance between the real roots inside $I_n$ of $P_n$ and $P_\infty$.  Via geometric considerations and properties of $h$,  it will follow that  $|N_n -N_\infty|=O(1)$ with high probability, which implies the desired estimates for $|N_n^*-N_\infty^*|$.

We'll use an elementary result whose proof is left as an exercise (see e.g. \cite{nnv2015}).
\begin{proposition}\label{p.rootcomp}
Let $F$   and $G$   be continuous real valued functions on $\R$, and $F\in C^2$. Let $\epsilon_1,M,N >0$ and $I:=[x_0-\epsilon_1/M, x_0+\epsilon_1/M]$. Assume that
\begin{itemize} 
\item $F(x_0)=0$, $|F'(x_0)| \ge \epsilon_1$;
\item  $|F''(x)| \le M$ for  $x\in I$;
\item  $\sup_{x\in I} |F(x)-G(x)| \le M'$.
\end{itemize}
Then $G$ has a root in $I$ if $M' \le \frac 1 4 \epsilon_1^2/M$.
\end{proposition}



\proof[Proof of Lemma~\ref{l.compareN}]

Let $h_t(x):=h(x/t)$, for $t>0$. 

Let $q_n(x)=e^{-x^2/2}P_n(x)$ and $q_\infty(x)=e^{-x^2/2}P_\infty(x)$. Note that the real roots of $q_n$ and $P_n$ are the same, and the real roots of $q_\infty$ and $P_\infty$ are the same.

Let $c_2>0$ and $c_1>c_2+2$. Let $I_n=[-n^{1/2}+n^{1/6}\log n, n^{1/2} - n^{1/6}\log n]$ and let $J_n =   \supp(h_{R_n}) \setminus I_n$.

Applying Lemma~\ref{l.derivatives-PnPinfty} (with $N=C_0\log^{1/2} n$, $C_0$ large), Lemma~\ref{l.drPn}, Lemma~\ref{l.drPinfty},
 with probablity $1-O(n^{-c_2})$  the following event (denoted by $E$) holds:  For every $x\in I_n$, we have

(i) $|q_n(x)-q_\infty(x)| \le M':= n^{-C}$

(ii) if $q_n(x)=0$ then $|q'_n(x)| \ge  \epsilon_1:=n^{-c_1}$ and $q_n(x')\ne 0$ for all $x'\in I_n$ such that   $|x-x'| \le \epsilon_1$.

(iii)   $|q_n''(x)|    \le M:= C_1 n \log ^{1/2}n$, $C_1$ absolute constant.

By choosing $C>2c_1+1$, it follows that
\begin{eqnarray*} 
\frac 1 4 \frac{\epsilon_1^2 }M  &=& \frac {n^{-2c_1-1} } {4C_1  \log n}  \\
&>& M' =n^{-C}
\end{eqnarray*}
Consequently,  Proposition~\ref{p.rootcomp} applies. (Note that the zeros of $P_n$ are at least $\epsilon_1$ apart by  (ii)). Thus for each zero of $P_n$ in $I_n$ (except for those near the endpoints) we could pair with one real zero of $P_\infty$ that is within a distance $\epsilon_1/M<\epsilon_1/2$.

Similarly, we consider the event $E'$ with $P(E') \ge 1-O(n^{-c_2})$ where  the following holds: for every $x\in I_n$,

(i) $|q_n(x)-q_\infty(x)| \le n^{-C}$.

(ii) If $q_\infty(x) = 0$ then $|q_\infty'(x)| \ge n^{-c_1}$ and $q'_n(x') \ne 0$ for all $x'\in I_n$ such that $|x-x'| \le e^{-c_1n}$.

(iii) $|q''_\infty(x)| \lesssim n\log^{1/2} n$.

Thus by applying Proposition~\ref{p.rootcomp} as before it follows that on the event $E'$ for each zero of $P_\infty$ in $I_n$ (except for those near the endpoints) we could pair with one real zero of $P_n$ that is within a distance $\epsilon_1/M<\epsilon_1/2$. 

Consequently, on the event $G=E' \cap E$ the zeros of $P_n$ and $P_\infty$ inside $I_n$ will form pairs, except for $O(1)$ zeros near the endpoints.

Now,  if $|x-x'|\le \epsilon_1/M$ is such a pair then there are three possibilities: (i) both $x$ and $x'$ are inside one interval forming $supp(h)$, or (ii) both $x$ and $x'$ are outside $supp(h)$, or (iii) one of them is inside and one is outside.

In the last two cases we have $|h(x/R_n)-h(x'/R_n)|=O(1)$, while in the first  case using H\"older continuity of $h$ we have 
$$|h(x/R_n)-h(x'/R_n)|\lesssim (\frac{\epsilon_1/M}{R_n})^\alpha < \frac 1 n$$ 
by choosing $c_2,c_1$ large compare to $1/\alpha$. Since there are at most $n$ such pairs, it follows that on the event $G$ we have
$$|N_\infty- N_n| \lesssim 1 + M_n + M_\infty \ \ ,$$
where $M_n$ and $M_\infty$ are the numbers of zeroes of $P_n$ and $P_\infty$ in $J_n$, respectively. 

Note that if $R_n\le n^{1/2}-n^{1/6}\log n$ then $supp( h_{R_n})\subset I_n$ and $M_n=M_\infty = 0$, so clearly $\E M_n^2 = \E M_\infty^2 = 0  = o(R_n)$. 
On the other hand, if $R_n > n^{1/2}-n^{1/6}\log n$ (recall that $R_n \le n^{1/2} + o(n^{1/4})$ by given assumption) then we have $|J_n| ^2 =o(n^{1/2}) = o(R_n)$. By translation invariance of the real zeros of $P_\infty$ and using Theorem~\ref{t.clt}, it follows that $$\E M_\infty^2 = O(|J_n|^2) = o(R_n)$$ 
By Lemma~\ref{l.localPn} we also have $M_n \le n^{o(1)}(1+|J_n \cap [-n^{1/2}, n^{1/2}]|) = O(n^{o(1)+1/6})$ with overwhelming probability. Since $M_n\le n$ always, it follows that
$$\E M_n^2 =   O(n^{1/3+o(1)}) = o(R_n)$$ 
Therefore,  taking $c_2$ large we obtain
\begin{eqnarray*} 
\E |N_n -N_\infty|^2 &\lesssim& 1 + \E M_\infty^2 + \E M_n^2 + \E((N_n^2+N_\infty^2)1_{G^c})\\
&\lesssim& o(R_n) + P(G^c)^{1/2} (\E[N_n^4] + \E [N_\infty^4])^{1/2}\\
&=& o(R_n) + O(n^{-c_2/2})O(n^4) = o(R_n) 
\end{eqnarray*}
here we have used the crude estimate $N_n= O(n)$ and the estimate  $\E N_\infty^4 \lesssim R_n^{4} =O(n^{2})$ (which is a result of Theorem~\ref{t.clt}).

It follows that in both cases we have
$$\E |N_n -N_\infty|^2  = o(R_n)$$
 
This completes the proof of Lemma~\ref{l.compareN}. \endproof

\section{Estimates for correlation functions}\label{s.correlation-functions}
In this section we summarize several new estimates  for the correlation function for the real zeros of $P_\infty$, which will be used in the proof of   Theorem~\ref{t.clt}. 

We first recall the notion of  correlation function. Let $X$ be a random point process on $\R$. For $k\ge 1$, the function $\rho:\R^k \to \R$ is the $k$-point corelation function of $X$ if for any compactly supported $C^\infty$ function $f:\mathbb R^k \to \R$ it holds that
$$\E \sum_{x_1,\dots, x_k} f(x_1,\dots, x_k) = \int\dots\int_{\R^k} f(\xi_1,\dots, \xi_k)\rho(\xi_1,\dots, \xi_k)d\xi_1\dots d\xi_k$$
where on the left hand side the summation is over all ordered $k$-tuples  of different elements in  $X$.\footnote{So if $(x_\alpha)_{\alpha\in I}$ is a labeling of elements of $X$ then we are summing over all $(x_{\alpha_1},\dots,x_{\alpha_k})$ where $(\alpha_1,\dots,\alpha_k)\in I^k$ such that $\alpha_i\ne \alpha_j$ if $i\ne j$. The correlation function is symmetric and the definition does not depend on the choice of the labeling.}   Note that this implies $\rho$ is locally integrable on $\R^k$. If there is $\epsilon>0$ such that $\rho$ is locally $L^{1+\epsilon}$ integrable, then by a simple approximation argument  it follows that   the above equality holds when   $f$ is only bounded and compactly supported. In particular, for every interval $I\subset \R$ it holds that
$$\E X_I(X_I-1)\dots (X_I-k+1) = \int\dots\int_{I^k} \rho(\xi_1,\dots,\xi_k) d\xi_1\dots d\xi_k$$
here $X_I = |X\cap I|$.

One should point out that the $k$-point correlation function does not always exists (however existence of the correlation measure, generalizing $\rho(\xi_1,\dots,\xi_k)d\xi_1\dots d\xi_k$, follows from the Riesz representation theorem). For the setting of the current work (namely $P_n$ and $P_\infty$), existence of the correlation functions is a consequence of the Kac-Rice formula (see also \cite{hkpv2009} for generalizations to all real Gaussian analytic functions).


Let $\rho^{(k)}$ be the $k$-point correlation function for the real zeros of $P_\infty$.

When it is clear from the context we will simply write $\rho$ instead of $\rho^{(k)}$.

\begin{lemma}\label{l.cor-local} For every $M>0$ and $k\ge 1$ there is a finite positive constant $C_{M,k}$ such that  for all $x_1,\dots, x_k \in  [-M,M]$ it holds that
$$\frac 1 {C_{M,k}} \prod_{1\le i<j\le k} |x_i-x_j| \le  \rho(x_1,\dots, x_k)  \le C_{M,k} \prod_{1\le i<j\le k} |x_i-x_j|$$
\end{lemma}

Lemma~\ref{l.cor-local} is a special case of the following more general  result, which holds for any $2k$-nondegenerate real Gaussian analytic functions on $\C$,  examples include random series $\sum_{j=0} c_j \xi_j x^j$ where $\xi_j$ are iid normalized Gaussian, $c_0, c_1,\dots \in \R$ such that $\sum_j c_j^2<\infty$ and $c_0,\dots, c_{2k-1}\ne 0$. This notion of nondegeneracy is a real variant of the complex nondegeneracy notion in \cite{ns2012}, see Section~\ref{s.background} of the current paper for details. 

\begin{lemma}\label{l.cor-gaf} Let $k\ge 1$. Let $f$ be a $2k$-nondegenerate real Gaussian analytic function on $\C$.  Let $\rho_f$ denote its $k$-point correlation function for the real zeroes. For every $M>0$  there is a finite positive constant $C_{M,k,f}$ such that  for all $x_1,\dots, x_k \in  [-M,M]$ it holds that
$$\frac 1 {C_{M,k,f}} \prod_{1\le i<j\le k} |x_i-x_j| \le  \rho_f(x_1,\dots, x_k)  \le C_{M,k,f} \prod_{1\le i<j\le k} |x_i-x_j|$$
\end{lemma}

Our next estimates will be about clustering properties for $\rho$.
\begin{lemma} \label{l.clustering}
There are finite positive constants $\Delta_k$ and $C_k$ such that the following holds: Given any $X=(x_1,\dots, x_k)$ of distinct points in $\R$, for any partition $X=X_I \cup X_J$ with $d=d(X_I,X_J)\ge 2\Delta_k$ we have
\begin{equation}\label{e.clustering}
|\frac{\rho(X)}{\rho(X_I)\rho(X_J)}-1| \le C_k \exp^{-\frac1 2 (d-\Delta_k)^2}
\end{equation}
\end{lemma}
Using  Lemma~\ref{l.cor-local}, it follows that if $X=(x_1,\dots, x_k)$ splits into two clusters $X_I$ and $X_J$ that are  sufficiently far part, then the correlation function essentially factors out.  From these clustering estimates and the well-known translation invariant properties of the real zeros of $P_\infty$ (see Lemma~\ref{l.Zinv} in Appendix~\ref{s.translation-invariant} for a proof), it follows that  $\rho$ is bounded globally. 

\begin{lemma} \label{l.cor-uniform} Let $\ell(t) = \min(1,|t|)$ for every $t\in \R$.  For every $k\ge 1$ there is a finite positive constant $C_k$ such that 
$$ \frac1{C_k} \prod_{1\le i<j\le k} \ell(|x_i-x_j|) \le \rho(x_1,\dots, x_k) \le C_k \prod_{1\le i<j\le k} \ell(|x_i-x_j|)$$
\end{lemma}
Indeed, if $k=1$ then the estimates hold trivially. The proof of the general case uses induction: for $k\ge 2$, if we could split  $X=(x_1,\dots, x_k$)  into two groups $X_1$, $X_2$ with distance $C\Delta_k$ where $C$ is sufficiently large (depending on $k$) then using Lemma~\ref{l.clustering} we have
$$|\frac{\rho(X)}{\rho(X_1)\rho(X_2)}-1| \le C_k e^{-\frac 1 2 (C-1)^2 \Delta_k^2} < \frac 1 2$$
therefore the desired claim follows from the induction hypothesis. If no such spliting could be found then  it follows from geometry that $diam(X)$ is bounded. Consequently, the desired bounds follow from the local estimates for the correlation function (Lemma~\ref{l.cor-local}) and the translation invariant properties of the real zeros.

Using Lemma~\ref{l.cor-uniform} and Lemma~\ref{l.clustering}, we immediately obtain the additive form of  \eqref{e.clustering}:
\begin{lemma}\label{l.clustering-add}
There are finite positive constants $\Delta_k$ and $C_k$ such that the following holds: Given any $X=(x_1,\dots, x_k) \in \R^k$, for any partition $X=X_I \cup X_J$ with $d=d(X_I,X_J)\ge 2\Delta_k$ we have
\begin{equation}\label{e.clustering-add}
|\rho(X) - \rho(X_I)\rho(X_J) | \le C_k \exp^{-\frac1 2 (d-\Delta_k)^2}
\end{equation}
\end{lemma}



\section{Proof of Theorem~\ref{t.clt} using correlation function estimates}\label{s.proof-clt}


Recall that $Z$ denotes the (multi-set of the) zeros of $P_\infty$ and  $h:\R \to \R_+$  is bounded and compactly supported, and
$$n(R,h) = \sum_{x\in Z} h(x/R)$$
for each $R>0$. Note that $n(R,h)=n(1,h_R)$ where $h_R(x)=h(x/R)$. For convenience of notation, let $\sigma(R,h)^2$ be the variance of $n(R,h)$ and let $n^\ast(R;h)$ be the normalization of $n(R,h)$, namely
$$n^\ast(R;h) = \frac{n(R,h)- \E n(R,h)}{\sigma(R,h)}$$

\subsection{Bound on the moments}\label{s.moment}

In this section, we'll show that $\E [n(R,h)^k] \lesssim R^k$. Cleary it suffices to consider $h=1_I$ for some fixed interval $I$.  Let $X$ denote $n(R,h)$ and let $I_R$ denote $\{Rx: x\in I\}$.   Using the uniform bound for the correlation function of real zeros of $P_\infty$ proved in Lemma~\ref{l.cor-uniform}, we have  
$$\E X(X-1)\dots (X-k+1)= \iint_{I_R\times \dots \times I_R} \rho(x_1,\dots, x_k)dx_1\dots dx_k$$
$$ \lesssim |I_R|^k (\min (1,diam(I_R))^{k(k-1)/2} \lesssim n^{k/2}$$
then the claims folow from writing $X^k$ as a linear combination of $X(X-1)\dots(X-j)$ with $j=0,1,\dots, k-1$.

\subsection{Asymptotic normality}

The convergence of $n^*(R,h)$ to standard Gaussian follows from the following two lemmas:
\begin{lemma}\label{l.clt-conditional}
Let $h:\R \to \R_+$ be bounded and compactly supported. Assume that there are $C,\epsilon>0$ such that $\sigma(R,h)\ge C R^{\epsilon}$ for $R$ sufficiently large. Then $n^\ast(R;h)$ converges in distribution to the standard Gaussian law as $R\to\infty$.  
\end{lemma}
\begin{lemma} \label{l.vargrow} If $h\in L^2$ then $\sigma(R,h)\gtrsim R^{1/2}$.
\end{lemma}
We will prove Lemma~\ref{l.vargrow} in Section~\ref{s.vargrow}. In this section we'll prove Lemma~\ref{l.clt-conditional}.

We will use the cumulant convergence theorem which will be recalled below.  The cumulants $s_k(N)$ of the random variable $N = \sum_{x\in Z} h(x)$ is defined by the formal equation
$$\log \E (e^{\lambda N}) = \sum_{k\ge 1} \frac {s_k(N)}{k!} \lambda^k$$
in particular
$$s_1(N) = \frac{d}{d\lambda}\Big( \log \E(e^{\lambda N})\Big)|_{\lambda=0} = \E N$$ 
$$s_2(N) = \frac{d^2}{d\lambda^2}\Big( \log \E(e^{\lambda N})\Big)|_{\lambda=0} = \E N^2 - (EN)^2 = \Var(N)$$

The version of the cumulant convergence theorem that we use is the following result of S. Janson \cite{janson1988}:

\begin{theorem}[Janson]\label{t.cct} Let $m\ge 3$. Let $X_1, X_2, \dots$ be a sequence of random variables such that  as $n\to\infty$ it holds that
\begin{itemize}
\item $s_1(X_n) \to 0$, and
\item $s_2(X_n) \to 1$, and
\item $s_j(X_n)\to 0$ for each $j\ge m$. 
\end{itemize}
Then $X_n\to N(0,1)$ in distribution as $n\to\infty$, furthermore all moments of $X_n$ converges to the corresponding moments of $N(0,1)$.
\end{theorem}


Since $s_1(N)  \equiv \E N$ and $s_2(N) \equiv  \Var(N)$ and  $n^*(R;h)$ has mean $0$ and variance $1$, it remains to show that the higher cumulants of $n^*(R;h)$ converge to $0$ as $R\to\infty$. We'll show that 
\begin{lemma} \label{l.sk-bound} For some finite constant $C_k$ depending only on $k$ it holds that 
$$s_k(n(1,h)) \le C_k \|h\|_\infty^k |supp(h)|$$ here $|supp(h)|$ is the Lebesgue measure of the support of $h$. 
\end{lemma}
Applying Lemma~\ref{l.sk-bound} to $h_R(x)=h(x/R)$, it follows that $s_k(n(R,h))   \le C_{h,k} R$. It follows from scaling symmetries and the definition of cumulants that if $N'=aN+b$ where $a>0$ and $b\in \R$ are fixed constants, then $s_k(N')=a^k s_k(N)$ for any $k\ge 2$. Thus, $s_k(n^*(R,h)) = \sigma(R,h)^{-j} s_k(n(R,h))$ for all $k\ge 2$. Consequently, using the fact that $\sigma(R,h)$ grows as a positive power of $R$, it follows that for $k$ sufficiently large $s_k(n^\ast(R;h)) \to 0$ as $R\to \infty$, as desired. 

Thus it remains to prove Lemma~\ref{l.sk-bound}. The proof uses the notion  of the truncated correlation functions , whose definition is recalled below. First, given  $Z=(x_1,\dots, x_k)$ let $|Z|:=k$ and let $Z_I$ denote $(x_j)_{j\in I}$. Let $\Pi(k)$ be the set of all partitions of $\{1,2,\dots,k\}$ (into nonempty disjoint subsets).  The truncated correlation function $\rho^T$ is defined using the following recursive formula (see e.g. Mehta \cite[Appendix A.7]{mehta1991}):
\begin{equation}\label{e.inverse-rho-rhoT}
\rho(Z) = \sum_{\gamma \in \Pi(k)} \rho^T(Z,\gamma)
\end{equation}
here $\Pi(k)$ is the set of partitions of $\{1,2,\dots,k\}$, and if $\gamma$ is the partition   $\{1,\dots,k\} = I_1\cup \dots \cup I_j$ then $\rho^T(Z,\gamma) := \rho^T(Z_{I_1}) \dots \rho^T(Z_{I_j})$. Explicit computation gives $\rho^T(x_1)=\rho(x_1)$, $\rho^T(x_1,x_2) = \rho(x_1,x_2) - \rho(x_1)\rho(x_2)$, etc.

To prove Lemma~\ref{l.sk-bound}, we will use the following two properties:

\begin{lemma}\label{l.sk-rhoT} For any $k\ge 1$ it holds that
\begin{equation}\label{e.cumulant}
s_k(n(1,h))= \sum_{\gamma \in \Pi(k)} \int_{\R^{|\gamma|}} h^\gamma(x)\rho^T(x)dx
\end{equation}
where $|\gamma|$ is the number of subsets in the partition $\gamma$,  $dA(x)$ is the Lebesgue measure on $\R^{|\gamma|}$, and if $\gamma_1,\dots, \gamma_j$ are the cardinality of the subsets in $\gamma$ then 
$$h^\gamma(x)=h(x_1)^{\gamma_1}\dots h(x_j)^{\gamma_j} $$
\end{lemma}

\begin{lemma} \label{l.rhoT-decay} There are finite positive constants $c_k, C_k$ such that for any $Z=(x_1,\dots,x_k)$ it holds that
$$\rho^T(Z) \le C_k \min(1, e^{-c_k |diam(Z)|^2})$$
\end{lemma}

A complex variant of Lemma~\ref{l.sk-rhoT} was also considered by Nazarov--Sodin in \cite{ns2012}, who provided a proof using a detailed algebraic computation.    Lemma~\ref{l.sk-rhoT} could be proved using a similar argument, and we include a proof in Appendix~\ref{s.sk-rhoT}.

\proof[Proof of Lemma~\ref{l.rhoT-decay}] We will use mathematical induction on $k$. If $k=2$ this follows from the uniform boundedness and clustering properties of $\rho$:
$$|\rho^T(x_1,x_2)|=|\rho(x_1,x_2) - \rho(x_1)\rho(x_2)| \le C\min(1, e^{-c|x_1-x_2|^2}) \ \ .$$

Let $k\ge 3$ and assume the estimates hold for all collection $k'$ points where $1\le k' <k$. Then there is a partition of $Z= Z_I\cup Z_J$ based on $\{1,\dots, k\}  = I \cup J$ such that $dist(Z_I,Z_J) \ge diam(Z)/Z_k$ and $I$ and $J$ are nonempty. It suffices to show that $|\rho^T(Z)| \le C_k e^{-c_k d(Z_I,Z_J)^2}$. 

Let $\Pi_1(k)$ be the set of partitions of $\{1,\dots,k\}$ that mixes $Z_I$ and $Z_J$, i.e. there is at least one block in the partition that intersects both $Z_I$ and $Z_J$, and let $\Pi_2(k)$ be $\Pi(k) \setminus \Pi_1(k)$. It follows from \eqref{e.inverse-rho-rhoT} that
$$\rho(Z) = \rho^T(Z) + \sum_{\gamma\in \Pi_1(k)}\rho^T(Z,\gamma) + \rho(Z_I)\rho(Z_J)$$
consequently using clustering of $\rho$ and the triangle inequality
$$|\rho^T(Z)| \le |\rho(Z) -  \rho(Z_I)\rho(Z_J)| + \sum_{\gamma\in \Pi_1(k)}|\rho^T(Z,\gamma)|$$
(and the induction hypothesis and boundedness of $\rho$)
$$\le C_k e^{-c_k dist(Z_I,Z_J)^2}$$
$$\le C_k e^{-c_k diam(Z)^2}$$
(note that the constants $c_k$ in different lines are not necessarily the same).
\endproof

We now finish the proof of Lemma~\ref{l.sk-bound}. Since $|\Pi(k)|\lesssim_k 1$, using Lemma~\ref{l.sk-rhoT} it suffices to show that
$$|\int_{\R^{|\gamma|}} h^\gamma(x)\rho^T(x)dA(x)| \le C_k \|h\|_\infty^k |supp(h)| $$
for each $\gamma\in \Pi(k)$. Fix such a $\gamma$.
Let $\gamma_1,\dots,\gamma_j$ be the length of its blocks. 
Using the uniform boundedness of  the correlation function $\rho^{(k)}$ (Lemma~\ref{l.cor-uniform}), we obtain
$$|\int_{\R^{|\gamma|}} h^\gamma(x)\rho^T(x)dA(x)| \le \int_{\R^j} |h(x_1)|^{\gamma_1}\dots |h(x_j)|^{\gamma_j} \rho^T(x_1,\dots, x_j)d x_1\dots dx_j$$
$$\le \|h\|_\infty^k |supp(h)| \sup_{x_1 \in \R} \int_{\R^{j-1}} |\rho^T(x_1,x_2,\dots,x_j)| dx_2\dots dx_j $$
$$\le C_k \|h\|_\infty^k |supp(h)| $$
in the last estimate we used Lemma~\ref{l.rhoT-decay}.

\endproof

\subsubsection{Growth of the variance}\label{s.vargrow}
In this section we prove Lemma~\ref{l.vargrow}. We have $\sigma(R,h) = \sigma(1,h_R)$, so we first estimate $\sigma(1,h)$ and then apply the estimate to $h_R$. Note that for $x\in \R$ we have $\rho(x) = \frac 1 \pi$. (See e.g. Edelman--Kostlan \cite{ek1995}). We then have
$$\sigma(1,h)^2=\iint_{\R^2} h(x_1)h(x_2)(\rho(x_1,x_2)-\rho(x_1)\rho(x_2) + \delta(x_1-x_2)\rho(x_1))dx_1dx_2$$
Let $k(x_1,x_2)=\rho(x_1,x_2)-\rho(x_1)\rho(x_2)$, since the distribution of the real zeros is translation invariant it follows that $\rho(x_1,x_2)$ depends only on $x_1-x_2$ (while $\rho(x_1)=\rho(x_2)=\frac 1 \pi$). Thus, we may write $k(x_1,x_2)=k(x_1-x_2)$ with
$$|k(x)| \le C e^{-Cx^2}$$
uniformly over $x\in \R$, thus in particular $\widehat k\in L^\infty\cap L^1$, and
$$\sigma(1,h)^2 = \iint_{\R^2} h(x_1)h(x_2)\Big[k(x_1-x_2)  + \frac 1 \pi \delta(x_1-x_2)\Big]dx_1dx_2$$
$$=\int_{\R} |\widehat h(\xi)|^2 [\frac 1\pi +\widehat k(\xi)] d\xi$$
Consequently,
\begin{eqnarray*}
\sigma(R,h)^2  &=& \int_{\R} |R\widehat h(R\xi)|^2 [\frac 1\pi +\widehat k(\xi)] d\xi \\
&=& R\int_{\R} |\widehat h(u)|^2 [\frac 1\pi +\widehat k(\frac u R)] du
\end{eqnarray*}
Using $\widehat k\in L^\infty$ and the dominated convergence theorem, it follows that
$$\lim_{R\to\infty} \frac{\sigma(R,h)^2}{R} = \lim_{R\to\infty} \int |\widehat h(u)|^2 (\frac 1 \pi + \widehat k(\frac uR)) du = (\frac 1 \pi + \widehat k(0)) \|h\|_2^2$$

Explicit computation \cite{sm2008}  gives $\widehat k(0) + \frac 1 \pi =0.18198... >0$ (for the reader's convenience we include a self-contained derivation in   Appendix~\ref{s.manju}).  Consequently $\sigma(R,h)\gtrsim R^{1/2}$.  


\endproof

\section{Real Gaussian analytic functions and linear functionals} \label{s.background}

In this section  we discuss   real Gaussian analytic functions and linear functionals on $\C$. These notions are adaptations of analogous notion in \cite{ns2012} and will be used in the proof of the correlation function estimates of Section~\ref{s.correlation-functions}.

\subsection{Real Gaussian analytic functions}

We say that $g$ is a real Gaussian analytic function (real GAF)  if  
$$g(z)=\sum_{j=1}^\infty \xi_j g_j(z)$$ 
where $\xi_i$ are iid normalized real Gaussian, and $g_1,g_2,\dots$ are  analytic functions   on $\C$ such that $\sum_j |g_j|^2<\infty$ uniformly on any compact subset of $\C$. In particular, uniformly over any compact subset of $\C$ we have $\E |g(z)|^2 = \sum_j |g_j(z)|^2< \infty$.

\subsection{Linear functionals}

We say that $L$ is a linear functional   if for some $K\ge 1$ there are  $m_1,\dots, m_K \in \mathbb Z$ nonnegative and $z_1,\dots, z_K \in \C$ and $\gamma_1,\dots, \gamma_K\in \C$ such that for any real GAF $g$   it holds that
$$Lg = \sum_{j=1}^K \gamma_j g^{(m_j)} (z_j) \ \ .$$
Here we require $(m_j,z_j) \ne (m_h,z_h)$ for $j\ne h$. We loosely say that  $z_j$ are the poles of $L$ (technically speaking only the distinct elements of $\{z_j\}$ should be called the poles of $L$, although in this case one has to count multiplicity).

Since $\sum_j |g_j|^2<\infty$ uniformly on compact subsets of $\C$, by standard arguments it follows that almost surely $\sum_{n=1}^\infty \xi_n g_n(z)$ converges absolutely on compact subsets of $\C$ to an analytic function (for a proof see e.g. \cite{hkpv2009}). Writing
\begin{eqnarray*}
Lg 
&=& \sum_{n\ge 1} \xi_n L(g_n) \ \ 
\end{eqnarray*}
and using independence of $\xi_n$'s, it follows that $Lg=0$ a.\,s. iff $L(g_n)=0$ for all $n$.


\subsection{Rank of linear functions}Let $L$ be a linear functional  with poles $z_1,\dots, z_K$.

Let $G \subset \C$   be a bounded domain with simple smooth boundary $\gamma=\partial G$ such that the poles $z_j$ are inside the interior $G^o$. 
By Cauchy's theorem, if $g$ is analytic then
\begin{equation}
\label{e.Lkernel}
 Lg=\int_\gamma g(z) r^L(z)dz \ \ , \ \ r^L(z) = \frac 1 {2\pi i} \sum_{j=1}^K \frac {\gamma_j m_j!}{(z-z_j)^{m_j+1}} \ \ .
\end{equation}
Now, $r^L$ is a rational function vanishing at $\infty$, and will be refered to as the  kernel of $L$.  We define the rank of $L$   to be the degree  of the denominator in any irreducible form of $r^L$ (this notion of rank is well defined and is independent of $G$).

\subsection{Degenerate and nondegenerate GAFs}
 
We say that a real GAF $g$   is  \emph{$d$-degenerate}   if  there is a linear functional $L\ne 0$  of rank at most $d$ such that $Lg = 0$ almost surely.   If no such linear functional exists, we say that $g$ is \emph{$d$-nondegenerate}.

\subsection{Linear functional arises from the Kac-Rice formula}

We discuss  linear functionals used in the proof.
Let $f$ be a real Gaussian analytic function.  Using the Kac-Rice formula for correlation functions of the real zeros of $f$  (see e.g. \cite{hkpv2009}) asserts that: for any $(x_1,\dots, x_k)\in \R^k$, we have
\begin{equation}\label{e.kac-rice}
\rho_k(x_1,\dots, x_k) = \frac{1}{(2\pi)^{k} |\det(\Gamma)|^{1/2}} \int_{\mathbb R^k} |\eta_1\dots \eta_k| e^{-\frac 1 2 \< \Gamma^{-1} \eta,\eta \>}d\eta_1\dots d\eta_k
\end{equation}
 where 
$\Gamma$ is the covariance matrix of 
$(f(x_1),f'(x_1),\dots, f(x_k), f'(x_k))$, and 
$$\eta=(0,\eta_1,\dots, 0,\eta_k) \in \mathbb R^{2k}$$


Given $\gamma=(\alpha_1,\beta_1,\dots, \alpha_k, \beta_k)\in \R^{2k}$, via elementary computations we have
\begin{equation}\label{e.bilinear}\<\gamma,\Gamma \gamma\> = \<\Gamma \gamma, \gamma\> = \E |Lf|^2 \ \ , \ \ Lf := \sum_{j=1}^k \alpha_j f(x_j) + \beta_j f'(x_j) \ \ .
\end{equation}
One could also define local version of $L$, namely  for any $I\subset \R$ we  $L_I f = \sum_{i \in I} \alpha_i f(x_i) + \beta_i f'(x_i)$ is also a linear functional.  Certainly  $L$ and $L_I$ depend on $\gamma$, however we will supress  the notational dependence for brevity, and none of the implicit constants in our estimates will depend on $\gamma$. Note that in the Kac-Rice formula, $\gamma=(0,\eta_1,\dots, 0,\eta_k)$.

\subsection{Non-degeneracy of random series}

Consider random infinite series $f(z)=\sum_{j\ge 0} a_j \xi_j z^j$ such that $(\xi_j)$ are iid standard Gaussian, $a_j\in \C$ and $\sup_{z\in K} \sum_j |a_jz^j|^2<\infty$ for any compact $K$. We now show that for such series if $a_0,\dots, a_{d-1}\ne 0$ then $f$ is $d$-nondegenerate on $\C$.  (Certainly $f$ is a  real GAF.)

Assume towards a contradiction that $f$ is $d$-degenerate. Then there is a linear functional $L$   of rank at most $d$ such that $L (a_n z^n)=0$ for all $n\ge 0$. Since $a_n\ne 0$ for $0\le n\le d-1$, it follows that $L(z^n)=0$ for all $0\le n\le d-1$. Taking $\gamma= \{|z|=R\}$ for any $R>0$ sufficiently large so that the poles of $L$ are enclosed inside $\gamma$,  we get
$$0= L (z^n) = \int_\gamma z^n r^L(z)dz$$
for all $n\in \overline{0,d-1}$. Since rank of $L$ is at most $d$ there is some $m\in \{1,d\}$ and $C\ne 0$ such that $z^m r^L(z) = C(1+o(1))$ as $|z|\to \infty$ uniformly. Consequently, 
$$ \int_{|z|=R} z^{m-1}r^L(z)dz = \int_\gamma z^m r^L(z) \frac {dz}z \to 2\pi i C \ne 0$$ as $R\to\infty$ contradiction.

It follows from the above discussion that the infinite flat series $P_\infty$ is $2k$-nondegenerate, and the Gaussian Kac polynomial of degree $2k-1$ defined by $g_{2k-1}(x)=\xi_0+\xi_1x+\dots + \xi_{2k-1}x^{2k-1}$ is also $2k$-nondegenerate.

\subsection{Equivalence of linear functionals}

The following lemma is a real Gaussian adaptation of a result in \cite{ns2012}.

\begin{lemma} \label{l.nond} Assume that $f$ is $d$-nondegenerate real GAF. Let $K\subset \C$ be nonempty compact.  Let $G$ be a bounded domain such that $K\subset G^o$, and assume that $\gamma=\partial G$ is a simple rectifiable curve. 

Then for any $d\ge 1$   there is a finite positive constant $C=C(d,G, K, f)$ such that for every linear functional $L$ of rank at most $d$  with poles in $K$ we have
$$\frac 1{C} \max_{z\in\gamma} |r^L(z)|^2 \le \E |Lf|^2 \le C \max_{z\in\gamma} |r^L(z)|^2$$
\end{lemma}

\proof   The proof largely follows \cite{ns2012}, included here to keep the paper self-contained.

We first show the upper bound. Let $ds$ denote the arclength measure on $\gamma$, then using \eqref{e.Lkernel} and Cauchy-Schwarz we have
$$|Lf|^2 \le length(\gamma)  \max_{z\in \gamma} |r^L(z)|^2 \int_{\gamma} |f(z)|^2 ds \ \ ,$$
$$\E|Lf|^2 \le \Big(length(\gamma) \max_{z\in \gamma} |r^L(z)|^2 \Big) \int_\gamma \E|f|^2 \le C_{\gamma,f} \max_\gamma |r^L|^2 \ \ .$$
We now show the lower bound. Assume towards a contradiction that the lower bound does not hold, then there is a sequence $(L_n)_{n\ge 1}$ of linear functionals of rank at most $d$ (with poles in $K$) such that $\max_\gamma |r_n|=1$ but
$$\lim_{n\to \infty} \E |L_n f|^2 = 0$$
We write $r_n(z)=\frac{p_n(z)}{q_n(z)}$ where $p_n$ and $q_n$ are  polynomials, and by multiplying both the numerator and denominator of $r_n$ by common factors (of the form $(x-\alpha)$ with $\alpha\in K$) if necessary we may assume that  $deg(q_n)=d$ and $deg(p_n)\le d-1$ and $q_n$ is monic. Since the zeros of $q_n$ are in $K$, we have $\sup_{z\in \gamma}|q_n(z)|<C_{d,K}<\infty$ (uniformly over $n$), therefore using $\sup_{\gamma} |r_n(z)|\le 1$ we obtain $\sup_{z\in \gamma} |p_n(z)| <C_{d,K}$ uniformly over $n$. Therefore, by passing to a subsequence, we may assume that $(p_n)$ converges uniformly on $\gamma$ to $p$. By iteratively  passing to subsequences we may assume further that $p_n', p_n'', \dots, p_n^{(d)}$ converge uniformly to $p', p'',\dots, p^{(d)}$. Since $\deg(p_n)<d$, it follows that $p^{(d)}\equiv 0$ and consequently $p$ is a polynomial of degree at most $d-1$. 

Now,  the $d$ complex zeros of $q_n$ are in $K$, a compact set, therefore by passing to a subsequence we may assume that uniformly on $\gamma =\partial G$ we have $q_n \to q$, and $q$ is a  monic polynomial of degree $d$ with zeros in $K$.  

Using partial fractional decomposition of $r(z)=p(z)/q(z)$, we obtain a linear functional $L$ of rank at most $d$ with poles in $K$ such that $\max_{z\in\gamma} |r^L(z)- r^{L_n}(z)| \to 0$ when $n\to\infty$. Consequently using the upper bound (already shown above) we obtain
$$\E |L_n f- Lf|^2 = O(\max_{\gamma}|r^{L_n}-r^L|^2) =o(1)$$
Using $lim_{n\to\infty} \E|L_n f|^2 = o(1)$ it follows that $\E |Lf|^2 = 0$, hence $Lf  = 0$ almost surely. This violates the $d$-nondegeneracy of $f$.
\endproof

\section{Local estimates for correlation functions}
In this section we prove Lemma~\ref{l.cor-local} and Lemma~\ref{l.cor-gaf}. Using Lemma~\ref{l.nond} and the Kac-Rice formula \eqref{e.kac-rice}, we observe that if $f_1$ and $f_2$ are two $2k$-nondegenerate real Gaussian analytic functions and $\rho^{[1]}$ and $\rho^{[2]}$ are the corresponding $k$-point correlation functions for the real zeroes,  then there is a finite positive constant $C=C_{M,N, k,f_1,f_2}$ such that $$\frac 1 {C} \rho^{[2]}(y_1,\dots, y_k) \le \rho^{[1]}(y_1,\dots, y_k) \le C \rho^{[2]}(y_1,\dots, y_k)$$

Indeed, let $\Gamma_j$  be the   covariance matrix for $f_j(y_1), f'_j(y_1),\dots, f_j(y_k), f'_j(y_k))$, which are positive definite symmetric. Then by Lemma~\ref{l.nond}, it follows that $\det(\Gamma_1)$ and $\det(\Gamma_2)$ are comparable and $\<\Gamma_1^{-1}\eta,\eta\>$ and $\<\Gamma_2^{-1}\eta,\eta\>$ are comparable. Consequently, using the Kac-Rice formula \eqref{e.kac-rice} it follows that $\rho^{[1]}(y_1,\dots, y_k) $ and $\rho^{[2]}(y_1,\dots, y_k)$ are comparable.

Therefore it suffices to show Theorem~\ref{l.cor-local}. Namely, we'll show that the correlation function for the real zeroes of $P_\infty$ is locally comparable to the Vandermonde product.

Let $M>0$ and $k\ge 1$. Assume that $x_1,\dots, x_k\in [-M,M]$. Let $N=N(M,k)$ be a large positive constant that will be chosen later. Thanks to the translation invariant property of the distribution of real zeros of $Z$, we have $\rho(x_1,\dots, x_k) = \rho(x_1+N,\dots, x_k+N)$. Let $y_1=x_1+N$, ..., $y_k=x_k+N$. Then $N-M\le y_j\le N+M$, and our choice of $N$ will ensure in particular that $N-M$ and $N+M$ are very large.

We now apply the above observation to $f_1=P_\infty$ and $f_2 = g_{2k-1}:=\xi_0+\xi_1 x+\dots + \xi_{2k-1}x^{2k-1}$, the Gaussian Kac polynomial. It then suffices to show that for any $n\ge k$ the correlation function $\rho_{Kac}$ for the real zeros of the Gaussian Kac polynomial $g_n(x)=\xi_0+\xi_1x+\dots + \xi_n x^n$ satisfies
$$\rho_{Kac}(y_1,\dots, y_k) \sim_{M,N,k,n} \prod_{1\le i\le j} |y_i-y_j|$$
whenever  $y_1,\dots,y_k \in [N-M,N+M]$ and $N-M\gg 1$.

We now observe that the distribution of the real roots of the Kac polynomial $g_n$ is invariant under the transformation  $x\mapsto 1/x$, underwhich $g_n(x)$ becomes $x^{-n} \widetilde g_n(x)$ where $\widetilde g_n=\xi_n+\xi_{n-1}x+\dots + \xi_0x^n$. 

It follows that, with $w_j=1/y_j$,
$$\rho_{Kac}(y_1,\dots, y_k) \sim_{M,N,k,n} \rho_{Kac}(w_1,\dots, w_k)$$
Indeed, it is well known (see e.g. \cite{hkpv2009}) that
\begin{eqnarray*} 
&\rho_{Kac}(y_1,\dots, y_k) & \\
&=& \lim_{\epsilon\to 0} \frac{\P(|u_1-y_1|\le \epsilon,\dots, |u_k-y_k|\le \epsilon \Big| g_n(u_1)=\dots=g_n(u_k)=0)}{(2\epsilon)^k}
\end{eqnarray*}
Now,  observe that if $|u_1-y_1|\le \epsilon$ and $\epsilon>0$ is sufficiently small then $|\frac 1 {u_1}-\frac 1 {y_1}| \lesssim \epsilon$ where the implicit constant depends on $M,N$. Conversely if $|\frac 1 {u_1}-\frac 1 {y_1}| \le \epsilon/C$ for $C$ very large depending on $M,N$ then for $\epsilon>0$ sufficienlty small we will have $|u_1-y_1|\le \epsilon$. It follows that 
$\rho_{Kac}(y_1,\dots, y_k)$ is comparable to the limit
$$\lim_{\epsilon \to 0} \frac{\P(|\frac 1{u_1}-w_1|\le \epsilon,\dots, |\frac 1{u_k}-w_k|\le \epsilon \Big| \widetilde g_n(\frac 1{u_1})=\dots=\widetilde g_n(\frac 1 {u_k})=0)}{(2\epsilon)^k}$$
which is exactly $\rho_{Kac}(w_1,\dots, w_k)$. 

Now, note that we also have
$$\prod_{1\le i<j\le k} |w_i-w_j| \sim_{M,N,k,n} \prod_{1\le i<j\le k} |y_i-y_j|$$
and note that $|w_j| \le \frac 1 {N-M}$ which could be made small if $N$ is chosen large. Therefore it suffices to show that for $\delta>0$ sufficiently small depending on $k$ and $n$ there is a finite positive constant $C=C_{\delta,k,n}$ such that for any $w_1,\dots, w_k \in [-\delta,\delta]$ it holds that
$$\frac 1 C \le \frac{\rho_{Kac}(w_1,\dots, w_k)}{\prod_{1\le i<j\le k}|w_i-w_j|} \le C$$
To show this estimate, our starting point is an explicit formula due to Gotze--Kaliada--Zaporozhets \cite[Theorem 2.3]{gkz2015}  for the real correlation of the general random polynomial
$$f(x)=\gamma_0+\gamma_1 x+\dots + \gamma_nx^n$$
where $\gamma_j$'s are independent real-valued random variables, and the distribution of $\gamma_j$ has probability density $f_j$. To formulate the formula, we first fix some notations. Given $w=(w_1,\dots, w_k)$ and $0\le i \le k$ we define $\sigma_i(w)$ to be the $i$th symmetric function of $x$, namely the sum of all products of $i$ coordinates of $w$: $$\sigma_i(w) =\sum_{1\le j_1<j_2<\dots<j_i\le k} w_{j_1}\dots w_{j_i}$$ 
(if $i>k$ or $k<0$ then $\sigma_i :=0$). Then we have, using \cite[Theorem 2.3]{gkz2015},

$$\rho(w_1,\dots,w_k)=\prod_{1\le i<j\le k} |w_i-w_j| \times$$
$$\times \int_{\mathbb R^{n-k+1}} \prod_{i=0}^n f_i  \Big(\sum_{j=0}^{n-k} (-1)^{k-i+j} \sigma_{k-i+j}(w)t_j\Big) \prod_{i=1}^k \Big|\sum_{j=0}^{n-k} t_j w_i^j\Big| dt_0 \dots dt_{n-k}$$
We apply this to $f=g_n$ the Gaussian Kac polynomial of degre $n$, note that $f_j(t)=\frac 1 {\sqrt{2\pi}}e^{-t^2/2} \lesssim 1$.

Note that if $\max |w_i|\le \delta$ for $\delta$ very small then for $k\le i \le n$ we have
$$\sum_{j=0}^{n-k} (-1)^{k-i+j} \sigma_{k-i+j}(w)t_j = t_{i-k} + O(\delta \max_j |t_j|) $$
therefore with $\delta$ small enough (depending on $k$ and $n$) 
$$\sum_{i=k}^n |\sum_{j=0}^{n-k} (-1)^{k-i+j} \sigma_{k-i+j}(x)t_j |^2  \sim_{n,k} \sum_{j=0}^{n-k} t_j^2$$
Since $0\le |\sum_{j=0}^{n-k} (-1)^{k-i+j} \sigma_{k-i+j}(x)t_j |^2 \lesssim_{k,n} \sum_{j=0}^{n-k} t_j^2$ for any $i$ (in particular for those $0\le i<k$), it follows immediately that
$$\sum_{i=0}^n |\sum_{j=0}^k (-1)^{k-i+j} \sigma_{k-i+j}(x)t_j |^2  \sim_{n,k} \sum_{j=0}^{n-k} t_j^2$$

Therefore for some finite positive constants $C_1,C_2$ that may depend on $k,n$ it holds that 
$$e^{- C_1 \sum_{j=0}^{n-k}t_j^2} \lesssim_{k,n} \prod_{i=0}^n f_i \Big(\sum_{j=0}^{n-k} (-1)^{k-i+j} \sigma_{k-i+j}(x)t_j\Big)  \lesssim_{k,n} e^{- C_2 \sum_{j=0}^{n-k}t_j^2}$$
From here it follows easily that
$$\rho(x_1,\dots, x_k) \lesssim_{\delta,k,n} \prod_{1\le i<j\le k} |x_i-x_j|$$

For the lower bound, note that if $t_1,\dots, t_{n-k} \in [-1,1]$ then $\prod_{i=1}^k |\sum_{j=0}^{n-k} t_j w_i^j| = |t_0|^k +O_{k,n}(\delta)$, so 
$$\int_{\R^{n-k+1}}\prod_{i=0}^n f_i  \Big(\sum_{j=0}^{n-k} (-1)^{k-i+j} \sigma_{k-i+j}(w)t_j\Big) \prod_{i=1}^k \Big|\sum_{j=0}^{n-k} t_j w_i^j\Big| dt_0 \dots dt_{n-k}$$
$$ \gtrsim_{k,n}  \int_{\R \times [-1,1]^{n-k}} (|t_0|^k  +O(\delta))e^{-C_1( t_0^2+\dots+t_{n-k}^2)}dt_0\dots dt_{n-k}$$
$$\gtrsim \int_{\R} |t_0|^k e^{-C_1 t_0^2} dt_0 + O(\delta) \gtrsim_{k,n} 1$$
if $\delta$ is sufficiently small.

This completes the proof of Lemma~\ref{l.cor-local}. \endproof

\section{Clustering estimates for correlation functions}
For a set $X=\{x_1,\dots, x_k\}$ of $k$ distinct points and for any nonempty subset $I\subset \{1,\dots, k\}$, we denote by $X_I$ the corresponding subset $\{x_i: i\in I\}$. Recall that $\rho$ denote the correlation function of the real zeroes of $P_\infty$. For simplicity of notation in this section let $f=P_\infty$.

In this section we will prove Lemma~\ref{l.clustering}, namely we will show that  there is a constant $\Delta_k>0$ and $C_k$ finite such that the following holds: Given any $X=(x_1,\dots, x_k)$ of distinct points in $\R$, for any partition $X=X_I \cup X_J$ with $d=d(X_I, X_J)\ge 2\Delta_k$ we have
$$|\frac{\rho(X)}{\rho(X_I)\rho(X_J)}-1| \le C_k \exp^{-\frac1 2 (d-\Delta_k)^2}$$

We will need the following lemma. Below fix $\eta=(0,\eta_1,\dots,0,\eta_k)$ where $\eta_1,\dots, \eta_k\in \R$, none of the implicit constants will depend on $\eta$. Let the linear functionals $L$ be defined using 
\begin{equation}\label{e.Ldef}
L f = \sum_{1\le j \le k} \eta_j f'(x_j)
\end{equation}
and define $L_I f$ for any $I\subset\{1,\dots, k\}$ using the summation over $j\in I$ instead of $1\le j\le k$.

\begin{lemma}\label{l.cov} There are finite positive constants $\Delta_k$ and  $C_k$  such that the following holds: Given any $X=(x_1,\dots, x_k)$ of distinct points in $\R$, for any partition $X=X_I \cup X_J$ with $d=d(X_I,X_J)\ge 2\Delta_k$ we have
$$|\E (L_I f) (L_J f)| \le C_k e^{-\frac 12 (d-\Delta_k)^2} (\E |L_I f|^2 + \E |L_J f|^2)$$
\end{lemma}
We defer the proof of this lemma to later sections. Below we prove the clustering property of the correlation function using this lemma.

Let $C_k$ and $\Delta_k$ be as in Lemma~\ref{l.cov}. 
Let $\epsilon = \frac 1 2 C_k e^{-\frac 1 2 (d-\Delta_k)^2}$ where $d=d(X_I, X_J)$. To show clustering it suffices to show that
\begin{equation}\label{e.mult-clustering}
(\frac{1-\epsilon}{1+\epsilon})^{k} \rho(X_I)\rho(X_J) \le \rho(X) \le (\frac{1+\epsilon}{1-\epsilon})^{k} \rho(X_I)\rho(X_J)
\end{equation}

Define $L$ using \eqref{e.Ldef}. By Lemma~\ref{l.cov} we have
$$|\E (L_I f) (L_J f)| \le \frac 1 2 \epsilon (\E |L_I f|^2 + \E |L_J f|^2)$$
therefore
$$(1- \epsilon) (\E |L_I f|^2 + \E |L_J f|^2) \le \E |L f|^2 \le (1+\epsilon) (\E |L_I f|^2 + \E |L_J f|^2)$$
Consequently using \eqref{e.bilinear} we obtain
\begin{equation}\label{e.gamma-comp}
(1-\epsilon) \Gamma_{I,J} \le \Gamma \le (1+\epsilon) \Gamma_{I,J}
\end{equation}
where $\Gamma$ is the covariance matrix of $(f(x_1),f'(x_1),\dots, f(x_k), f'(x_k))$ and $\Gamma_{I,J} = \begin{pmatrix} \Gamma_I & 0\\ 0 & \Gamma_J\end{pmatrix}$. 

We obtain $\det(\Gamma) \ge (1-\epsilon)^{2k} \det(\Gamma_I)\det(\Gamma_J)$, and therefore
$$\rho(X) \le (1-\epsilon)^{-k} \frac 1 { (2\pi)^k|\det(\Gamma_{I,J})|^{1/2} } \int_{\R^k} |\eta_1\dots \eta_k| e^{-\frac 1 2 (1+\epsilon)^{-1}\<\Gamma_{I,J}^{-1}\eta,\eta\>}d\eta$$
$$=(1-\epsilon)^{-k} (1+\epsilon)^{k} \frac 1 {(2\pi)^{|I|} |\det(\Gamma_{I})|^{1/2}} \int_{\R^{|I|}} |\eta_1\dots \eta_{|I|}| e^{-\frac 1 2 (1+\epsilon)^{-1}\<\Gamma_{I}^{-1}\eta,\eta\>}d\eta$$
$$\times
 \frac 1 {(2\pi)^{|J|}  |\det(\Gamma_{J})|^{1/2}} \int_{\R^{|J|}} |\eta_1\dots \eta_{|J|}| e^{-\frac 1 2 (1+\epsilon)^{-1}\<\Gamma_{J}^{-1}\eta,\eta\>}d\eta$$
$$=(\frac{1+\epsilon}{1-\epsilon})^{k} \rho(X_I)\rho(X_J)$$
Similarly we have
$$\rho(X) \ge (\frac{1-\epsilon}{1+\epsilon})^{k} \rho(X_I)\rho(X_J) \ \ .$$
This completes the proof of \eqref{e.mult-clustering}. 
\endproof

\subsection{Proof of Lemma~\ref{l.cov}}


We first show a small scale version of the lemma. 
\begin{lemma} Let $\rho>0$. Suppose that $K_1$ and $K_2$ are two intervals with length at most $2\rho$. Assume that $L_{K_j}$ is a linear functional on $\C$ with poles inside $K_j$ with  rank at most $k$. Assume that $d=dist(K_1, K_2) \ge 2\rho$. Then
$$\E(L_{K_1} f)  (L_{K_2}f) \le C_{f,k,\rho} e^{-\frac 1 2 (d-2\rho)^2} (\E |L_{K_1} f|^2 + \E |L_{K_2} f|^2)$$
\end{lemma}

\proof Let $c_1$ and $c_2$ be the centers of $K_1$ and $K_2$. 
Let $T_x f$ denote 
$$T_x f(z) = f(z+x) e^{-(\frac 1 2 |x|^2 + xz)}$$
then $\E T_xf(t_1)  T_xf(t_2) = \E f(t_1) f(t_2)$ for any $t_1,t_2 \in \R$. Therefore for any $x\in \R$, $T_x f$ and $f$ have the same distribution (in particular the distribution of the real zeroes of $f$ is translation invariant). Let $K_3$ and $K_4$ be $K_1 - c_1$ and $K_4=K_2 -c_2$, thus these intervals are centered at $0$ and have length at most $2\rho$. Let $L_1$ be such that $L_1 T_{c_1} f = L_{K_1}f$ and  let $L_2$ be such that $L_2 T_{c_2}f = L_{K_2}$. Then it is clear that $L_1$ and $L_2$ are linear functional of rank at most $k$ with poles inside $[-\rho,\rho]$. Let $\gamma =\{|z|=2\rho\}$. Since $f$ is $2k$-nondegenerate on $\C$, we then have
$$\E (L_{K_1} f L_{K_2}f) =  \E (L_{1} T_{c_1}f) (L_{2}T_{c_2}f)$$
$$\le |\int_{\gamma} \int_{\gamma} r^{L_1}(z_1) r^{L_2}(z_2) \E (T_{c_1}f)(z_1)T_{c_2}f(z_2))dz_1dz_2|$$
$$\le C_{\rho} \sup_{z_1,z_2 \in \gamma} |\E (T_{c_1}f (z_1)T_{c_2}f(z_2))|  (\max_{\gamma} |r^{L_1}|^2 + \max_{\gamma} |r^{L_2}|^2)$$
By explicit computation and using $|c_1-c_2|\ge 2\rho$, we obtain
$$\sup_{z_1,z_2\in\gamma}\E (T_{c_1} f(z_1)T_{c_2}f(z_2)) \le C_\rho e^{-(|c_1-c_2| - 2\rho)^2} $$
therefore
$$\E (L_{K_1} f L_{K_2}f) \le  C_\rho e^{-(|c_1-c_2| - 2\rho)^2} (\max_{\gamma} |r^{L_1}|^2 + \max_{\gamma} |r^{L_2}|^2) $$
$$\le C_{f,k,\rho}  e^{-|d - 2\rho|^2} (\E |L_1 f|^2 + \E |L_2 f|^2)$$
Since $f$ and $T_{c_j}f$  have the same distribution, the last display is the same as
$$= C_{f,k,\rho}  e^{-|d - 2\rho|^2} (\E |L_{K_1}f|^2 + \E |L_{K_2} f|^2)$$
which implies the desired estimate.
\endproof

We now start the proof of Lemma~\ref{l.cov}. We will construct a covering $X$ by small intervals having the following properties:

(i) The cover will consists of $m\le k$ intervals $I_1$,\dots $I_m$ each of length at most $\rho$  such that the distance between the centers of any two of them is at least  $4\rho$.   

(ii) The algorithm will ensure that $\rho>1$ (or any given large absolute constant) but $\rho = O_{k,f}(1)$. \\

We first let $\rho_1=1$  and use the given points as centers of the interval. 

If there are two centers with distance not larger than $4\rho_1$,  we replace these two centers by one center at their midpoint, and enlarge all intervals, replacing $\rho_1$ by $\rho_2=3\rho_1$. 

We repeat this process if needed, and since there are only $k$ points the process has to stop. Clearly the last radius is at most $3^{k-1}$. 

Note that we could ensure that $\rho$ is larger than any given constant depending on $k,f$ if needed, by setting $\rho_1$ the initial radius to be larger than this constant.

Now, choose $\Delta_k>4\rho$ such that $\Delta_k-4\rho$ is very large compared to $1$.

Notice that there is no $k$ such that $I_k$ intersects both $I$ and $J$. Let $A=\{k: I_k\cap I \ne \emptyset\}$ and $B=\{k: I_k \cap J \ne \emptyset\}$. Now using the above small scale result we have

$$\E |L_{I}f|^2 = \sum_{k\in A} \E |L_{I_k} f|^2 + \sum_{k, n \in A: k\ne n} \E (L_{I_k} f L_{I_n}f) \ge \frac 1 2 \sum_{k\in A} \E |L_{I_k} f|^2 $$
Similarly, $\E |L_{J}f|^2  \ge \frac 1 2 \sum_{k\in B} \E |L_{I_k} f|^2$.

Now, if $k\in A$ and $n\in B$ it is clear that $dist(I_k, I_n) \ge d-2\rho > 2\rho$ since $d\ge 2\Delta_k$ is very large compared to $\rho$. Therefore
$$|\E (L_I f)(L_Jf)| \le C_{f,k,\rho}\sum_{k\in A, n\in B} e^{-\frac 1 2 |d-4\rho|^2} (\E|L_{I_k} f|^2 + \E |L_{I_n} f|^2)$$
$$\le C_{f,k,\rho} e^{- \frac 1 2 |d-4\rho|^2} (\E |L_{I}f|^2  + \E |L_{J}f|^2 )$$
\endproof

\appendix

\section{Translation invariance of the real zeros of $P_\infty$}\label{s.translation-invariant}
The following property is standard, we include a proof for the convenience of the reader.
\begin{lemma}\label{l.Zinv} The distribution of the real roots of $P_\infty$ is invariant under translations on $\R$ and the reflection $x\mapsto -x$.
\end{lemma}
\proof Notice that $f(x)$ is a real Gaussian process with correlation function
$$K_f(x,y) = \E f(x) f(y) = \sum_{k \ge 0} \frac{x^k y^k}{k!} = e^{xy}$$
Let $g(x)=e^{-bx + \frac 1 2 b^2} f(ax+b)$ where $a\in \{-1,1\}$ and $b\in \R$. Then $g$ is also a centered real Gaussian  process with the correlation function
$$K_g(x,y) = K_f(ax+b, ay+b) = e^{(ax+b)(ay+b) - b(ax+b) - b(ay+b) + b^2}$$
$$=e^{xy}= K_f(x,y)$$
It follows that $g$ has the same distribution as $f(x)$. Consequently the real zeros of $f(ax+b)$ has the same distribution as the real zeros of $f(x)$.
\endproof

\section{Relation between cumulants and truncated correlation functions}\label{s.sk-rhoT}

For the convenience of the reader, we include a self-contained proof of  Lemma~\ref{l.sk-rhoT} in this section, which largely follows an argument in Nazarov--Sodin \cite{ns2012}. Recall that $X$ is the random point process for the real zeros of $P_\infty$, $h:\R\to \R_+$ is bounded compactly supported, and $n(1,h)=\sum_{\alpha\in X}h(\alpha)$, and $\Pi(k)$ is the collection of all partition of $\{1,\dots, k\}$ into disjoint nonempty subsets, and for each $\gamma \in \Pi(k)$ let  $|\gamma|$ be the number of subsets in the partition $\gamma$, and if $\gamma_1,\dots, \gamma_j$ are the cardinality of the subsets in  $\gamma$ then 
$$h^\gamma(x)=h(x_1)^{\gamma_1}\dots h(x_j)^{\gamma_j} $$

We first prove an analoguous relation between the moment and the (standard) correlation functions: if $m_k(N)=\E N^k$ denotes the $k$th moment of the random variable $N$, then
\begin{equation}\label{e.mk-rho}
m_k(n(1,h)) = \sum_{\gamma\in \Pi(k)} \int_{R^{|\gamma|}} h^{\gamma}(x) \rho(x)dx  \  \ .
\end{equation}
Indeed, divide $\Pi(k)$ into $\displaystyle \bigcup_{j\ge 1} \Pi(k,j)$ basing on $|\gamma|=j$. For each $\gamma\in \Pi(k,j)$ let $A_1,\dots, A_j$ be the subsets in the partition and let $\gamma_j=|A_j|$. Note that the correlation functions are uniformly bounded therefore we may use the bounded compactly supported function $H(x_1,\dots,x_j)=h(x_1)^{\gamma_1}\dots h(x_j)^{\gamma_j}$ as a test function in the defining property of correlation functions, and obtain
\begin{eqnarray*}
\int_{\R^{|\gamma|}} h^\gamma(x)dx 
&=& \int_{\R^j} h(x_1)^{\gamma_1}\dots h(x_j)^{\gamma_j}\rho(x_1,\dots,x_j)dx_1\dots dx_j \\
&=& \E \sum_{(\xi_1,\dots,\xi_j)} h(\xi_1)^{\gamma_1}\dots h(\xi_j)^{\gamma_j}
\end{eqnarray*}
where the last summation is over all ordered $j$-tuples of different elements of $X$. By summing over all possible values of $j$ and $\gamma_1,\dots,\gamma_j\ge 1$ (with $\gamma_1+\dots+\gamma_j=k$), it follows that
\begin{eqnarray*}
\sum_j \sum_{\gamma\in \Pi(k,j)}  \int_{\R^{|\gamma|}} h^\gamma(x)dx 
&=& \E \sum_{j}\sum_{\gamma_1+\dots+\gamma_j=k}\sum_{(\xi_1,\dots,\xi_j)} h(\xi_1)^{\gamma_1}\dots h(\xi_j)^{\gamma_j}\\
&=& \E [\sum_{\xi\in X} h(\xi)]^k \ \ =  \ \  m_k(n(1,h)) \ \ .
\end{eqnarray*}

We now use induction on $k$ to prove \eqref{e.cumulant}. Clearly \eqref{e.cumulant} holds for $k=1$. The key ingredient for the induction step is the following relationship between cummulant and moments (see e.g. \cite[Chapter 2]{shiryaev1996})
$$s_k = m_k - \sum_{j\ge 2} \sum_{\pi \in \Pi(k,j)} s_{\pi_1}\dots s_{\pi_j}$$
which is analogous to the following reformulation of \eqref{e.inverse-rho-rhoT}:
$$\rho^T(Z) = \rho(Z) - \sum_{j\ge 2} \sum_{\gamma\in \Pi(k,j)} \rho^T(Z,\gamma)$$
where $\rho^T(Z,\gamma) = \rho^T(Z_{\Gamma_1}) \dots \rho^T(Z_{\Gamma_j})$ if $\gamma=(\Gamma_1,\dots,\Gamma_j)$.

To facilitate the notation, for $\gamma, \pi\in \Pi(k)$ we say that $\gamma \le \pi$ if $\gamma$ is a refinement of $\pi$, in other words the partitioning subsets of $\gamma$ are subsets of the partitioning subsets in $\pi$.  If $\gamma\le \pi$ and $\gamma \ne \pi$ we say $\gamma<\pi$.

 Let $1$ denote the trivial partition  with just one partitioning subset, clearly all $\pi\in\Pi(k)$ satisfies $\pi \ll 1$. We'll write $\pi=(\Pi_1,\dots,\Pi_{|\pi})$ below.

Using the induction hypothesis we have
\begin{eqnarray*}
s_k 
&=& m_k - \sum_{\pi < 1} \prod_{j=1}^{|\pi|} s_{\pi_j}\\
&=& \sum_{\gamma  \le 1} \int_{\R^{|\gamma|}} h^\gamma(x)\rho (x)dx - \sum_{\pi < 1} \prod_{j=1}^{|\pi|} \sum_{\Gamma(j) \text{ partitions }  \Pi_j}\int_{\R^{|\Gamma(j)|}} h^{\Gamma(j)}(x^{(j)})\rho^T(x^{(j)})dx^{(j)} \ \ .
\end{eqnarray*}
Note that $x^{(j)}$ is a vector in $\R^{|\Gamma(j)|}$. Interchanging the sum in the second term and let $\gamma=(\Gamma(1),\dots,\Gamma(|\pi|))\le \pi$, we obtain
\begin{eqnarray*}
s_k 
&=&  \sum_{\gamma  \le 1} \int_{\R^{|\gamma|}} h^\gamma(x)\rho(x)dx  -  \sum_{\gamma  \le 1} \sum_{\gamma\le \pi <1} \int_{\R^{|\gamma|}} h^{\gamma}(x) \Big[\prod_{j=1}^{|\pi|} \rho^T(x^{(j)})dx^{(j)}  \Big]  \\
&=& \sum_{\gamma\le 1} \int_{\R^{|\gamma|}} h^{\gamma}(x) \Big[\rho(x)-\sum_{\pi: \gamma\le \pi <1} \prod_{j=1}^{|\pi|}  \rho^T(x^{(j)})dx^{(j)}\Big] \ \ , \ \ \text{here $x=(x^{(1)},\dots, x^{(|\pi|)})$,} \\
&=& \sum_{\gamma\le 1} \int_{\R^{|\gamma|}} h^{\gamma}(x) \Big[\rho(x)-\sum_{\pi\in \Pi(|\gamma|): \ \pi <1}  \rho^T(x, \pi)\Big]dx \\
&=& \sum_{\gamma\le 1} \int_{\R^{|\gamma|}} h^{\gamma}(x) \rho^T(x) dx \ \ .
\end{eqnarray*}
This completes the induction step and   the proof of Lemma~\ref{l.sk-rhoT}.\endproof

\section{An explicit computation for the two-point correlation function of $P_\infty$}\label{s.manju}

First thanks to translation and reflection invariance one has $\rho(s,t)=\rho(0,t-s) = \rho(0,s-t)$, therefore it suffices to compute $\rho(0,t)$ for $t>0$. This function was computed   in an earlier paper \cite{sm2008} (which also contains many interesting statistics about the real roots); we choose   provide 
the details for the reader's convenience.

Let $\gamma(t)=e^{-t^2/2}$ and $g(t)=\gamma(t) P_\infty(t)$. The zero distribution of $P_\infty$ and $g$ are the same, so it suffices to compute the two-point correlation function for the real zeros of $g$.  The covariance matrix for $g(0), g(t), g'(0), g'(t)$ is a symmetric $4\times 4$ matrix $\begin{pmatrix}A & B\\ C&D\end{pmatrix}$ where $A,B, C,D$ are $2\times 2$ matrices. It follows that the conditional distribution of $(g'(0), g'(t))$ given $g(0)=0$ and $g(t)=0$ is a centered bivariate Gaussian with covariance matrix $\Sigma = D - C A^{-1}  B$. Since $\E [g(t)g(s)] = \gamma(t-s)$, one has $C=B^T$ and
$$A= \begin{pmatrix}1 & \gamma(t) \\ \gamma(t) & 1\end{pmatrix} \ \ , \ \ B =   \begin{pmatrix}0 & -t\gamma(t) \\ t\gamma(t) & 0\end{pmatrix} \ \ ,  \ \ D =  \begin{pmatrix}1 & (1-t^2)\gamma(t) \\ (1-t^2)\gamma(t) & 1\end{pmatrix}\ \ .$$
Therefore via explicit computation (below $\gamma=\gamma(t) = e^{-t^2/2}$)
\begin{equation}\label{e.xydef}
\Sigma = I+M \ \ , \ \  M  = \begin{pmatrix} x& y\\y & x\end{pmatrix} \ \ , \ \ x = \frac{-t^2 \gamma^2}{1-\gamma^2} \ \ ,  \ \ y = (1-t^2)\gamma - \frac{t^2 \gamma^3}{1-\gamma^2}
\end{equation}
Note that $(g(0),g(t))$ is a centered bivariate Gaussian with covariant $A$, thus the density at $(0,0)$ is $(2\pi \sqrt{1-\gamma^2})^{-1}$. It follows that
\begin{eqnarray}\label{e.rho0t}
\rho(0,t) = \frac{1}{2\pi  \sqrt{1-\gamma^2}}\E[|X| \cdot |Y|]
\end{eqnarray}
where $(X,Y)$ have mean zero bivariate Gaussian distributions with covariant $\Sigma$.

Note that $|\alpha| = \frac 1 {\sqrt{2\pi}} \int_0^\infty (1-e^{-\alpha^2 x/2}) x^{-3/2}dx$ (which follows from  a simple rescaling of the integration variable and integration by parts). Using the Kac-Rice formula, it follows that
$$\E[|X|\cdot |Y|] =  \frac 1{2\pi} \int_0^\infty \int_0^\infty (f(u,v)-f(u,0)-f(0,v)+f(0,0)) u^{-3/2}v^{-3/2}dudv$$
$$f(u,v)  = \E [e^{-(uX^2+vY^2)/2}] = \Big|\det(I+\Sigma \begin{pmatrix}u & 0 \\ 0 & v\end{pmatrix}\Big|^{-1/2} $$
$$=\Big((1+u)(1+v) + x(u+v+2uv) + uv(x^2-y^2)\Big)^{-1/2}$$
Redefine $A=x$ and $B=x^2-y^2$, we have 
\begin{eqnarray*}
f(u,v) 
&=& \Big(1+(A+1)v + (1+A+v+2Av+Bv)u\Big)^{-1/2}\\ 
&=&  \Big(1+(A+1)v\Big)^{-1/2} \Big(1+ M_v u\Big)^{-1/2} 
\end{eqnarray*}
where $M_v=\frac{1+A+(1+2A+B)v}{1+(A+1)v}$. We now use the elementary identity
$$\int (1+\alpha)^{-1/2}\alpha^{-3/2}d\alpha = -2(1+\alpha)^{1/2} \alpha^{-1/2}+C$$
and obtain, via  a simple change of variables,
\begin{eqnarray*}
\int_0^\infty \Big(f(u,v) - f(0,v)\Big)u^{-3/2} du 
&=& \frac 1{\sqrt{1+(A+1)v}}  \int_0^\infty \Big((1+M_v u)^{-1/2} - 1\Big)u^{-3/2} du\\
&=& \frac {1}{\sqrt{1+(A+1)v}} \Big(2u^{-1/2} - 2 (1+ M_v u)^{1/2}u^{-1/2}\Big)|_0^\infty\\
&=& \frac {-2M_v^{1/2} }{\sqrt{1+(A+1)v}}
\end{eqnarray*}
In particular we could let $v=0$ and obtain $-2(1+(A+1)u)^{-1/2}$ as the result. Therefore
\begin{eqnarray*}
\E[|X|\cdot |Y|] 
&=& \frac 1 {2\pi} \int_0^\infty 2\Big(\sqrt{A+1} - \frac{(1+A+(1+2A+B)v)^{1/2}}{1+(A+1)v} \Big) v^{-3/2}dv
\end{eqnarray*}
Let $N= \frac {1+2A+B}{(1+A)^2}<1$ (since $B<A^2$). Using the change of variables $v\mapsto (1+2A+B)v/(1+A)$ we obtain 
\begin{eqnarray*}
\E[|X|\cdot |Y|]  &=& \frac {1}{2\pi}  2(1+2A+B)^{1/2}  \int_0^\infty \Big(1- \frac{(1+ v)^{1/2}}{1+N^{-1}v} \Big) v^{-3/2}dv\\
&=& \frac {1}{2\pi}  2(1+2A+B)^{1/2}  \Big(\frac 2 \alpha \arctan(\frac u \alpha) +\frac 2 u - 2v^{-1/2}\Big)_{v=0}^{v=\infty}
\end{eqnarray*}
where $u:=v^{1/2}(1+v)^{-1/2}$ and $\alpha = \sqrt{N/(1-N)}$. Thus
$$\E[|X|\cdot |Y|]  = \frac 1 {2\pi} 4(1+2A+B)^{1/2} \Big(\frac 1 \alpha \arctan(\frac 1 \alpha) + 1\Big)$$
Note that $\arctan(\frac 1 \alpha) = \arccos(\alpha/\sqrt{1+\alpha^2}) = \arccos(\sqrt N)$. Therefore
$$\E[|X|\cdot |Y|]  = \frac 1 {2\pi} \Big(4(1+2A+B)^{1/2}  + 4\sqrt{A^2-B}\arccos(\frac{\sqrt{1+2A+B}}{1+A})\Big)$$
 Recalling $A=x$ and $B=x^2-y^2$ and \eqref{e.xydef} we obtain
\begin{eqnarray}\label{e.rho0t-final}
\rho(0,t)  = \frac 1{\pi^2 \sqrt{1-e^{-t^2}}}  \Big(\sqrt{(1+x)^2-y^2} + |y| \arcsin(\frac {|y|}{1+x})\Big)
\end{eqnarray}
where  $x,y$ are defined using \eqref{e.xydef}. One could check that $1+x\ge 0\ge y$ thus by letting 
$$\delta:=\frac{|y|}{1+x} = \frac{e^{-t^2/2}(e^{-t^2/2}+t^2-1)}{1-e^{-t^2}-t^2e^{-t^2}}$$ 
we obtain
$$\rho(0,t)=\frac{\sqrt{(1-e^{-t^2})^2-t^4 e^{-t^2}}}{\pi^2(1-e^{-t^2})}\Big(1+\frac{\delta}{\sqrt{1-\delta^2}}  \arcsin\delta\Big)$$
which recovers \cite[(D8,D9)]{sm2008}. From here a numerical evaluation gives

\begin{corollary}\label{c.2cor-explicit} Let $k(t) = \rho(0,t)-\frac 1 {\pi^2}$. Then $\widehat k(0)+\frac1 \pi = 0.18198...>0$.
\end{corollary}

\begin{figure}\label{l.kfigure}
\centering
\includegraphics[scale=1]{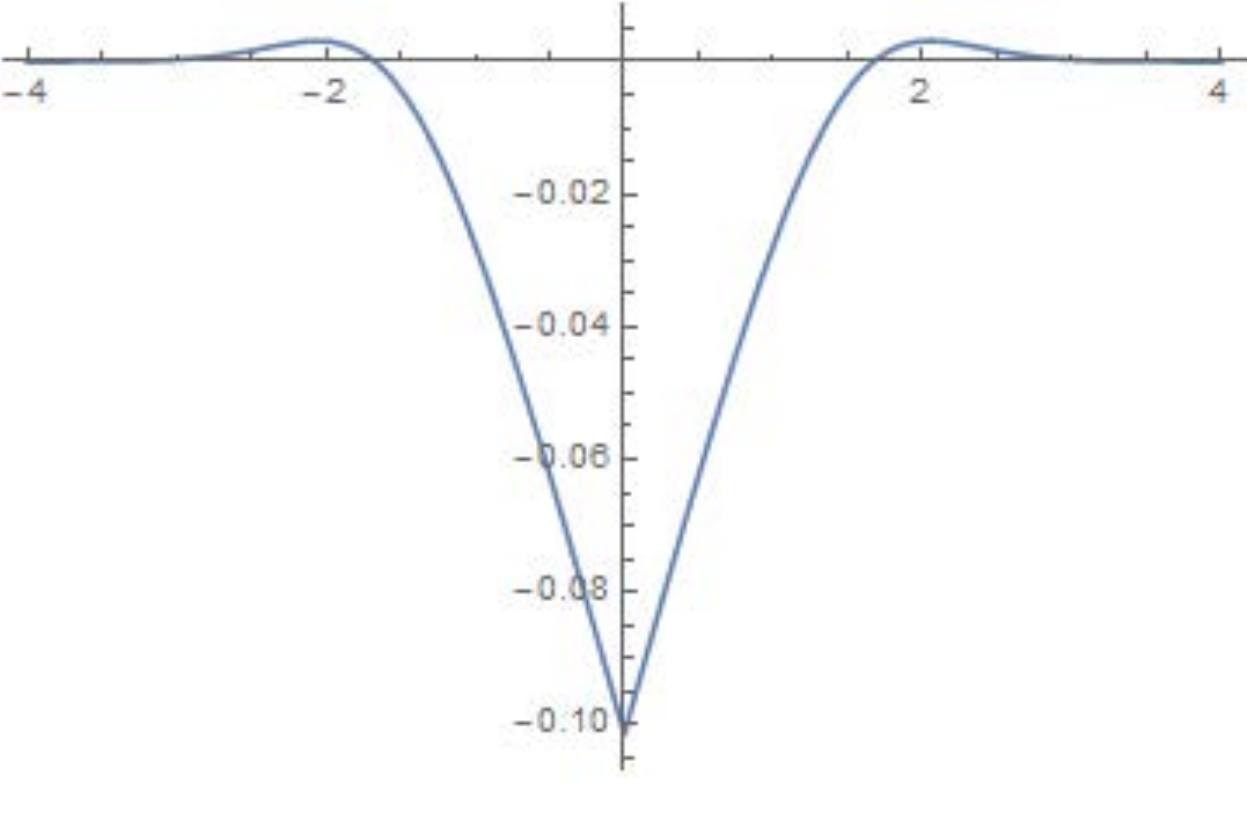}
\caption{Mathematica plot of $k(t):=\rho(0,t)-\frac 1 {\pi^2}$}
\end{figure}

 

\bibliographystyle{plain}
\end{document}